\documentclass{article}
\usepackage{tikz}
\usepackage[utf8]{inputenc}
\usepackage[T1]{fontenc}
\usepackage{mathtools}
\usepackage{hyperref}
\usepackage{stmaryrd}
\usepackage{amsmath}
\usepackage{amssymb}
\usepackage{mathrsfs}
\usepackage{graphicx}
\usepackage{MnSymbol,bbding,pifont}
\usepackage{yfonts}
\usepackage{amsthm}
\usepackage{tikz-cd}
\usepackage{authblk}
\usepackage{verbatim}
\usepackage[a4paper, total={6in, 8in}]{geometry}
\newtheorem{theorem}{Theorem}[section]
\newtheorem{lemma}[theorem]{Lemma}
\newtheorem{corollary}[theorem]{Corollary}

\newtheorem{proposition}[theorem]{Proposition}
\newtheorem{definition}[theorem]{Definition}
\newtheorem{example}[theorem]{Example}
\newtheorem{remark}[theorem]{Remark}

\newtheorem{problem}[theorem]{Problem}

\numberwithin{equation}{section}
\newenvironment{prf}{{\noindent\it Proof. }}{\hfill $\clubsuit $\par}

\newcommand\db{{\color{blue} \char`\\}}
\newcommand\nsb{{\color{red}/}}
\newcommand{\langlands}[1]{#1^{\mathbb{L}}}
\mathcode`*=\string"8000
\begingroup
\catcode`*=\active
\xdef*{\noexpand\textup{\string*}}
\endgroup
\title{Intersections of twisted cotangent bundles and \\ symplectic duality}
\author{Naichung Conan Leung and Yunsong Wei}
\date{}

\begin{document}
\setcounter{section}{-1}
\maketitle
\begin{abstract}
We observe that numerous symplectic resolutions can be expressed as intersections of twisted cotangent bundles. Additionally, their dual symplectic resolutions can be derived from intersections of dual twisted cotangent bundles. We determine the collection of fixed points for certain intersections that are Poisson slices, extending the computations of fixed points for parabolic Slodowy varieties.
\end{abstract}

\tableofcontents

\section{Introduction}Given a holomorphic symplectic manifold $X$, Rozansky-Witten \cite{rozansky1997hyper} defined a 3d B-model with target $X$. When $X$ is a symplectic resolution of symplectic singularity $\bar{X}$, we expect that
there is also a 3d A-model with target $X$. Furthermore, 3d mirror symmetry \cite{intriligator1996mirror},
also known as symplectic duality, predicts that there is another symplectic resolution ${X}^!\rightarrow \bar{X}^!$ which switches the A-model and B-model between $X$ and $X^!$. We will also refer to these as 3d mirror varieties, as these pairs are predicted to arise from physics as the mirror symmetry of the Higgs/Coulomb branches of a 3d \( \mathcal{N} = 4 \) supersymmetric quantum field theory \cite{braverman2016towards,kamnitzer2022symplectic,nakajima2015towards, nakajima2022mathematical}.

There are many surprising mathematical predictions arising from such a physical duality (see e.g. \cite{aganagic2021elliptic,botta2023bow, braden2014quantizations,braden2012quantizations, bullimore2017coulomb,bullimore2016boundaries, gammage2022perverse, hoang2024around}). The simplest one is a bijection between the fixed
point set of a torus action on $X$ and the corresponding one for $X^!$.   Many interesting algebraic varieties are studied within this class, such as parabolic Slodowy varieties \cite{hoang2024around} (also referred to as S3-varieties in \cite{braden2014quantizations,webster2011singular}), bow varieties \cite{ji2023bow,nakajima2017cherkis,rimanyi2020bow}, universal implosions \cite{dancer2024complex, fu2025affine,gannonproof,gannon2023differential,jia2022geometry}, and universal centralizers \cite{bezrukavnikov2005equivariant}. Notably, all these varieties are symplectic reductions of $T^*G$, the cotangent bundle of a reductive group \( G \), and their 3d mirror varieties can also be obtained as symplectic reductions of $T^*\langlands{G}$, where \( \langlands{G} \) is the Langlands dual group. 

In \cite{chan20243d}, the authors interpreted 3d mirror symmetry via 2d mirror symmetry for intersections of 3d branes. The basic philosophy is if \( \underline{C}_1 \) and \( \underline{C}_2 \) are two 3d branes in \( X \) with mirror 3d branes \( \underline{C}_1^! \) and \( \underline{C}_2^! \) in \( X^! \), then the intersections \( \underline{C}_1 \times_{X} \underline{C}_2 \) and \( \underline{C}_1^! \times_{X^!} \underline{C}_2^! \) should be 2d mirror to each other. Namely, their symplectic geometry and complex geometry got switched under mirror symmetry. In the same spirit, SYZ proposal \cite{strominger1996mirror} describes 2d mirror symmetry as T-duality, which is a duality in 1d TQFT. In this article, we aim to understand 3d mirror symmetry from the perspective of intersections of 4d branes in S-dual 4d gauge theories for \( G \) and \( \langlands{G} \). A 4d brane of \( G \) gauge theory in this context refers to certain holomorphic symplectic manifolds \( X \) endowed with a Hamiltonian \( G \)-action, which we denote as \( \underline{X} = (G \curvearrowright X \rightarrow \mathfrak{g}^*) \)\footnote{Certain anomaly cancellation condition of Gaiotto-Witten \cite{Gaiotto:2008ak} needed to be satisfied, which we will omit here.}. We primarily consider those which are twisted cotangent bundles \( T^*G / \!\!/_{\psi} H \) in this article.

S-duality predicts that given a 4d brane $\underline{X}$ for $G$ gauge theory, there
is another 4d brane, denoted as $\langlands{\underline{X}}$, for $\langlands{G}$ gauge theory. Furthermore, if $\underline{X}_1$ and $\underline{X}_2$ are two 4d branes, then
 \( \underline{X}_1 \underset{[\mathfrak{g}^*/G]}{\times} \underline{X}_2 \)  and $\langlands{\underline{X}}_1 \underset{[{\langlands{\mathfrak{g}}}^*/\langlands{G}]}{\times} \langlands{\underline{X}}_2$ should be 3d mirror to each other, where \( \underline{X}_1 \underset{[\mathfrak{g}^*/G]}{\times} \underline{X}_2 \) is the quotient of the fiber product $X_1\times_{\mathfrak{g}^*}X_2$ by $G$ as defined in Chapter \ref{Intersection of Hamiltonian $G$-varieties}.\footnote{If we view \( \underline{X}_i, i=1,2 \) as one-shifted Lagrangians \( [X_i/G] \rightarrow [\mathfrak{g}^*/G] \) in the one-shifted symplectic  stack \( [\mathfrak{g}^*/G] \), then \( \underline{X}_1 \underset{[\mathfrak{g}^*/G]}{\times} \underline{X}_2 \) could be viewed as the shifted Lagrangian intersection of \( [X_1/G] \) and \( [X_2/G] \). See \cite{safronov2017symplectic}.} This is closely related to the relative Geometric Langlands conjecture proposed and studied by Ben-Zvi, Sakellaridis and Venkatesh in \cite{ben2024relative}, where, in particular, they constructed a   candidate for \(\langlands{\underline{X}}\) when \(\underline{X}\) is a twisted cotangent bundle of a spherical variety. When \(\underline{X} = (G \curvearrowright T^{\ast} \mathbf{N} \rightarrow \mathfrak{g}^*)\) where \(\mathbf{N}\) is a smooth affine $G$-variety, in \cite{Nakajima:2024mlb}, Nakajima proposed a construction of \(\langlands{\underline{X}}\) generalizing the Braverman–Finkelberg–Nakajima (BFN) construction of Coulomb branches for 3d gauge theory with matters. We shall use \(\mathbf{M}\) and \(\langlands{\mathbf{M}}\) to refer to the pairs constructed by them in \cite{ben2024relative} and \cite{Nakajima:2024mlb}.

Fix a base $\Pi$ of the root system $\Phi$ of $G$. Given any subset $\mathcal{I}\subset\Pi$, let $P_\mathcal{I}\supset L_\mathcal{I}\supset B_{L_\mathcal{I}}$ be the standard parabolic subgroup associated to $\mathcal{I}$, the Levi subgroup of \( P_\mathcal{I} \) and the Borel subgroup of \( L_\mathcal{I} \) respectively. Let \( {U}_{P_\mathcal{I}} \) be the unipotent radical of \( P_\mathcal{I} \). In particular,  ${U}_{P_\mathcal{\varnothing}}=U_{B}=U$ is the unipotent radical of the Borel subgroup $B$. Let \( {e}_\mathcal{I} \) be the regular nilpotent element of \( L_\mathcal{I} \), and let \( {U}_{e_\mathcal{I}} \) be the corresponding unipotent group  defined in Section \ref{section bielawski slices}. For example, if \( G = \operatorname{GL}_4 \), then the base is \( \Pi = \{\alpha_1, \alpha_2, \alpha_3\} \) with \( \alpha_i = E_{i,i+1} \). Let \( \mathcal{I} = \{\alpha_1, \alpha_3\}\subset \Pi \); we have the following subgroups and the nilpotent element:
$$ 
\begin{array}{ccccccc}
     P_{\mathcal{I}} &\supset &L_{\mathcal{I}}  &\supset &B_{L_{\mathcal{I}} } & &U_{P_{\mathcal{I}}}  \\
    \begin{pmatrix}
         *& * &  *  & * \\  * & *&  * & * \\0&0&  * & * \\0&0&  * & * \\
     \end{pmatrix} & & \begin{pmatrix}
         *& * &  0  & 0\\  * & *& 0 & 0\\0&0&  * & * \\0&0&  * & * \\
     \end{pmatrix}& & \begin{pmatrix}
         *& * &  0 & 0 \\  0 & *&  0 & 0 \\0&0&  * & * \\0&0& 0 & * \\
     \end{pmatrix}& &\begin{pmatrix}
         1& * &  0 & 0 \\  0 & 1&  0 & 0 \\0&0&  1 & * \\0&0& 0 & 1 \\
     \end{pmatrix}\\
     &&\text{\rotatebox[]{90}{$\in$}}&&&&\\
     &e_{\mathcal{I}}=&\begin{pmatrix}
         0& 1 &  0 & 0 \\  0 & 0&  0 & 0 \\0&0& 0 & 1 \\0&0& 0 & 0 \\
     \end{pmatrix}&&&U_{e_\mathcal{I}}=&\begin{pmatrix}
         1& 0 &  0 & 0 \\ * & 1&  * & 0 \\0&0&  1 & 0 \\ * &0& * & 1 \\
     \end{pmatrix}
\end{array}
$$
\newpage
We suggest the following list of 4d mirror branes (also known as S-dual branes in \cite{ben2024relative}) in Chapter \ref{Section 4d mirror} as follows,  where the first row (i) is essentially suggested by Braden, Licata, Proudfoot and Webster in \cite{braden2014quantizations}\footnote{In \cite{braden2014quantizations}, the authors refer to the intersections of 4d branes in row $(i)$ as $S3$ varieties and investigate the 3d mirror pairs of $S3$ varieties. We reinvestigate them from the 4d perspective.} and the second row is also suggested by Tom Gannon and Sanath Devalapurkar in \cite[Remark~2.9]{gannon2025cotangent} and \cite{devalapurkar2024ku}.
\begin{equation}
\label{list}
\renewcommand{\arraystretch}{1.5}
 \begin{array}{cc|c}&\underline{X} & \langlands{\underline{X}}\\ 
\hline (i)&T^*G/\!\!/P_\mathcal{I}  & T^*\langlands{G}/\!\!/_{\langlands{e}_\mathcal{I}} \langlands{U}_{e_\mathcal{I}}  \\(ii)&
T^*G/\!\!/L_\mathcal{I} & T^*\langlands{G}/\!\!/_{\langlands{e}_\mathcal{I}} \langlands{U}\\(iii)&
T^*G/\!\!/B_{L_\mathcal{I}} & T^*\langlands{G}/\!\!/{U}_{\langlands{P}_\mathcal{I}} 
\end{array}   
\end{equation}
For the twisted cotangent bundles \( T^*G / \!\!/_{\psi} H \) listed above, the Hamiltonian action $G\curvearrowright  T^*G / \!\!/_{\psi} H\rightarrow \mathfrak{g}^*$ is implicitly understood, and so they are viewed as 4d branes.  The varieties listed in the second column should be interpreted as affinization of the corresponding twisted cotangent bundle when the context involves mirror symmetry. 

For the examples mentioned previously, we have: (a) the parabolic Slodowy varieties and bow varieties are the intersections of $T^*G/\!\!/P_\mathcal{I}$ and $T^*{G}/\!\!/_{{e}_\mathcal{J}} {U}_{e_\mathcal{J}}$ for $\mathcal{I},\mathcal{J}\subset \Pi$, whose suggested mirror 4d branes are $T^*\langlands{G}/\!\!/_{\langlands{e}_\mathcal{I}} \langlands{U}_{e_\mathcal{I}}$ and $T^*\langlands{G}/\!\!/\langlands{P}_\mathcal{J}$,
(b) the universal implosion is the affinization of the intersection of $T^*G/\!\!/U$ and $T^*G$, whose 4d mirror branes are suggested to be $T^*\langlands{G}/\!\!/{\langlands{T}}$ and $T^*\langlands{G}/\!\!/{\langlands{B}}$, where $\langlands{T}:=\langlands{L}_\varnothing$ is the maximal torus of $\langlands{G}$, 
(c) the universal centralizer is the self-intersection of $T^*{G}/\!\!/_{{e}_{\Pi}} {U}_{e_{\Pi}}$, whose 4d mirror brane is suggested to be $T^*\langlands{G}/\!\!/{\langlands{G}}=\operatorname{pt}$, a point.

Let us further elaborate on the relationships of these 4d branes listed in (\ref{list}) with 3d mirror symmetry in the following examples.

(1) By taking the intersection of the mirror 4d branes in row $(i)$ of the list (\ref{list}), the 3d mirror variety of the parabolic Slodowy variety 
\[
T^*G/\!\!/P_\mathcal{I} \underset{[\mathfrak{g}^*/G]}{\times} T^*{G}/\!\!/_{{e}_\mathcal{J}} {U}_{e_\mathcal{J}}
\]
is suggested to be another parabolic Slodowy variety
\[
T^*\langlands{G}/\!\!/_{\langlands{e}_\mathcal{I}} \langlands{U}_{e_\mathcal{I}} \underset{[{\langlands{\mathfrak{g}}}^*/\langlands{G}]}{\times} T^*\langlands{G}/\!\!/\langlands{P}_\mathcal{J}.
\]
In fact, a refined Nakajima-Hikita conjecture is proved for them by Hoang, Krylov, and Matvieievskyi in \cite{hoang2024around}. For type \( A \), the duality between their elliptic stable envelopes is proved by Botta and Rimanyi in \cite{botta2023bow}. When \(\mathcal{I} = \mathcal{J} = \varnothing\), we simply get the 3d mirror varieties $T^*G/\!\!/B$ and $T^*\langlands{G}/\!\!/\langlands{B}$.
The Koszul duality for Category $\mathcal{O}$ in type $A$ is proved in \cite{braden2014quantizations, kamnitzer2019category}. Other expected 3d mirror relationships between \( T^*G/\!\!/B \) and \( T^*\langlands{G}/\!\!/\langlands{B} \) are also described in \cite{dinkins20223d, rimanyi2019three}. 

(2) Also, row $(i)$ of list (\ref{list}) suggests that the 3d mirror variety of the intersection
$$T^*G/\!\!/P_\mathcal{I} \underset{[\mathfrak{g}^*/G]}{\times} T^*{G}/\!\!/{P}_{\mathcal{J}}$$ should be 
$$T^*\langlands{G}/\!\!/_{\langlands{e}_\mathcal{I}} \langlands{U}_{e_\mathcal{I}} \underset{[{\langlands{\mathfrak{g}}}^*/\langlands{G}]}{\times} T^*\langlands{G}/\!\!/_{\langlands{e}_\mathcal{J}} \langlands{U}_{e_\mathcal{J}}.$$ We have an explanation for this duality in type $A$ through bow diagrams and bow varieties in Section \ref{section bow variety}.

(3) As suggested by row $(ii)$ of this list (\ref{list}), the 3d mirror variety of 
\[
T^*G/\!\!/L_\mathcal{I} \underset{[\mathfrak{g}^*/G]}{\times} T^*{G}/\!\!/L_\mathcal{J}
\]
should be 
\[
T^*\langlands{G}/\!\!/_{\langlands{e}_\mathcal{I}} \langlands{U} \underset{[{\langlands{\mathfrak{g}}}^*/\langlands{G}]}{\times} T^*\langlands{G}/\!\!/_{\langlands{e}_\mathcal{J}} \langlands{U}.
\]
When \(\mathcal{I} = \mathcal{J} = \Pi\), this conjecture is supported by the fact that the Coulomb branch of the pure $G$ gauge theory, i.e. \( G \curvearrowright \{0\} \) is the universal centralizer of \(\langlands{G}\) \cite{bezrukavnikov2005equivariant, braverman2016towards, Nakajima:2024mlb}.

(4) We can also take the intersections of 4d branes from different rows in the list (\ref{list}). For example, this suggests that the 3d mirror variety of 
\[
T^*G/\!\!/P_\mathcal{I} \underset{[\mathfrak{g}^*/G]}{\times} T^*{G}/\!\!/_{{e}_\mathcal{J}} {U}
\]
is 
\[
T^*\langlands{G}/\!\!/_{\langlands{e}_\mathcal{I}} \langlands{U}_{e_\mathcal{I}} \underset{[{\langlands{\mathfrak{g}}}^*/\langlands{G}]}{\times} T^*\langlands{G}/\!\!/_{\langlands{e}_\mathcal{J}} \langlands{U}.
\]
When \(\mathcal{J} = \Pi\), these two are parabolic Slodowy varieties. When \(\mathcal{J} = \varnothing\), we calculated the set-theoretical fixed point set of 
\[
T^*G/\!\!/P_\mathcal{I} \underset{[\mathfrak{g}^*/G]}{\times} T^*{G}/\!\!/{U}
\]
in Proposition \ref{NegaDimension fixed proposition} and the fixed point set of 
\[
T^*\langlands{G}/\!\!/_{\langlands{e}_\mathcal{I}} \langlands{U}_{e_\mathcal{I}} \underset{[{\langlands{\mathfrak{g}}}^*/\langlands{G}]}{\times} T^*\langlands{G}/\!\!/{\langlands{T}}
\]
in Corollary \ref{Levi fixed point set}. There is a bijection between these two fixed point sets when \(\mathcal{I} = \varnothing\). Notice that $\mathcal{I}=\varnothing$ if and only if the dimension, $\dim G-\dim P_\mathcal{I}-\dim U$, of the quotient 
\(
T^*G/\!\!/P_\mathcal{I} \underset{[\mathfrak{g}^*/G]}{\times} T^*{G}/\!\!/{U}
\) is non-negative. The list also suggests that the 3d mirror variety of 
\[
T^*G/\!\!/P_\mathcal{I} \underset{[\mathfrak{g}^*/G]}{\times} T^*{G}/\!\!/L_{\mathcal{J}}
\]
is 
\[
T^*\langlands{G}/\!\!/_{\langlands{e}_\mathcal{I}} \langlands{U}_{e_\mathcal{I}} \underset{[{\langlands{\mathfrak{g}}}^*/\langlands{G}]}{\times} T^*\langlands{G}/\!\!/_{\langlands{e}_\mathcal{J}} \langlands{U}.
\]
When \(\mathcal{I} = \varnothing\) and \(G\) is of type \(A\), this is supported by the work of Tom Gannon and Harold Williams\cite{gannon2023differential} as will be explained in Section \ref{Section intersection quiver}.

(5) For \(\underline{X}\) in the list (\ref{list}), we also suggest that the intersection 
\[
\underline{X} \underset{[{{\mathfrak{g}}}^*/{G}]}{\times} \mathbf{M}
\]
is 3d mirror to 
\[
\langlands{\underline{X}} \underset{[{\langlands{\mathfrak{g}}}^*/\langlands{G}]}{\times} \langlands{\mathbf{M}},
\]
where \((\mathbf{M}, \langlands{\mathbf{M}})\), as mentioned before, is a pair of Hamiltonian dual hyperspherical varieties constructed by Ben-Zvi, Sakellaridis, and Venkatesh in \cite{ben2024relative} or S-dual pairs constructed by Nakajima in \cite{Nakajima:2024mlb}. In particular, when $\mathbf{M}=T^*\mathbf{N}$ and $\langlands{\mathbf{M}}$ is the S-dual constructed in \cite[Section~3.4]{Nakajima:2024mlb},
\[
T^*\mathbf{N} \underset{[{{\mathfrak{g}}}^*/{G}]}{\times} T^*G/\!\!/G = T^*\mathbf{N}/\!\!/G
\]
should have a 3d mirror variety:
\[
\langlands{\mathbf{M}} \underset{[{\langlands{\mathfrak{g}}}^*/\langlands{G}]}{\times} T^*\langlands{G}/\!\!/_{\langlands{e}_{\Pi}} \langlands{U} = \langlands{\mathbf{M}}/\!\!/_{\langlands{e}_{\Pi}} \langlands{U}
.\]
Sanath Devalapurkar provides a proof for this 3d mirror pair, which is presented in Proposition \ref{proposition Sympl dual}.

In the list (\ref{list}), \( T^*{G}/\!\!/L_{\mathcal{J}} \) is the cotangent bundle of a smooth affine variety \( {G}/L_{\mathcal{J}} \). We expect that the the affinization of the suggested 4d mirror brane 
\(
T^*\langlands{G}/\!\!/_{\langlands{e}_\mathcal{I}} \langlands{U}
\)
coincides with Nakajima's construction \( \langlands{\mathbf{M}} \) in \cite{Nakajima:2024mlb} for \( \mathbf{M} = T^*{G}/\!\!/L_{\mathcal{J}} \). In (\ref{list}), \( T^*G/\!\!/P_\mathcal{I} \) and \( T^*{G}/\!\!/_{{e}_\mathcal{I}} {U} \) are twisted cotangent bundles of spherical varieties. We expect that the 4d mirror branes to \( T^*G/\!\!/P_\mathcal{I} \) and \( T^*{G}/\!\!/_{{e}_\mathcal{I}} {U} \) as suggested in (\ref{list}) match with the corresponding constructions by BSV in \cite[Part~1, Section~4]{ben2024relative}.

Given nilpotent elements \( e,f \in \mathfrak{g} \), let \( \bar{\mathcal{O}}_e \) be the nilpotent orbit closure of \( e \) and let \( S_f \) be the Slodowy slice associated with \( f \). We have the identification 
\(
G \times S_f \cong T^*G/\!\!/_f U_f.
\)
We can also view \( \bar{\mathcal{O}}_e \) as a 4d brane \( (G \curvearrowright \bar{\mathcal{O}}_e \hookrightarrow \mathfrak{g}^*) \), and the nilpotent Slodowy slice \cite{arakawa2024singularities} is an intersection of 4d branes given by
\[
S_{e,f} := \bar{\mathcal{O}}_e \cap S_f = \bar{\mathcal{O}}_e \underset{[{{\mathfrak{g}}}^*/{G}]}{\times} (G \times S_f).
\]
Physics literature such as \cite{cabrera2019quiver,cabrera2017nilpotent,cremonesi2015t,hanany2020quiver} suggests that the 3d mirror to \( S_{e,f} \) is \( S_{\operatorname{d}_{BV}(f),\operatorname{d}_{BV}(e)} \), where \( \operatorname{d}_{BV}:\mathcal{N}/G \rightarrow \langlands{\mathcal{N}}/\langlands{G} \) is the Barbasch-Vogan map \cite{barbasch1985unipotent} that maps from the set of nilpotent orbits of \( \mathfrak{g} \) to the set of nilpotent orbits of \( \langlands{\mathfrak{g}} \). In type \( A \), the affinization of \( T^*G/\!\!/P_\mathcal{I} \underset{[\mathfrak{g}^*/G]}{\times} T^*{G}/\!\!/_{{e}_\mathcal{J}} {U}_{e_\mathcal{J}} \) is a nilpotent Slodowy slice. More generally, the affinization of \( T^*G/\!\!/P_\mathcal{I} \underset{[\mathfrak{g}^*/G]}{\times} T^*{G}/\!\!/_{{e}_\mathcal{J}} {U}_{e_\mathcal{J}} \) admits a finite map to a nilpotent Slodowy slice \cite{braden2014quantizations}. The Barbasch-Vogan map is compatible with the duality suggested in the first row of the list (\ref{list}) \cite[Section~5]{hoang2024around}.

As mentioned before, a particularly nice property of 3-dimensional mirror varieties is the existence of a bijection between the fixed-point sets under certain torus actions. To verify the conjectural 3d mirror varieties suggested above and to expand it, we attempted to find the 4d mirror branes corresponding to \( T^*G/\!\!/H \), where \( H \) is a subgroup with maximal rank in \( G \). For example, $T^*G/\!\!/P_\mathcal{I}, T^*G/\!\!/L_\mathcal{I}$, and $T^*G/\!\!/B_{L_\mathcal{I}}$ in the left column of the list (\ref{list}) belong to this class. Therefore, it motivates us to calculate the fixed-point set of \( T^*G/\!\!/H \underset{[\mathfrak{g}^*/G]}{\times} 
 T^*G/\!\!/_{{e}_\mathcal{I}} {U}_{e_\mathcal{I}} \), which are Poisson slices in \( T^*G/\!\!/H \) in the sense of \cite{crooks2020log}. We construct explicitly a torus action on $ T^*G/\!\!/H \underset{[\mathfrak{g}^*/G]}{\times} 
 T^*G/\!\!/_{{e}_\mathcal{I}} {U}_{e_\mathcal{I}}$ in Chapter \ref{Section Some fixed points subsets}. The torus is the product of the maximal torus of the centralizer of \( e_{\mathcal{I}} \) and another 1-dimensional torus \( \mathbb{C}^\times \). We show that such Poisson slices have discrete fixed points (Lemma \ref{general fixed point}) under this torus action.

A connected subgroup containing the maximal torus of $G$ corresponds to a closed subset \( \Gamma \) of the root system \( \Phi \). For example, the Levi subgroup $L_{\mathcal{I}}$ corresponds to $\Phi_\mathcal{I}$, the subset of $\Phi$ generated by elements in $\mathcal{I}$; the subgroup $B_{L_{\mathcal{I}}}$ corresponds to $\Phi_\mathcal{I}^+$, the set of positive roots in $\Phi_\mathcal{I}$. We denote by \( H_\Gamma \) the corresponding connected subgroup of $G$ associated to the closed subset \( \Gamma \) of $\Phi$.

Given $\mathcal{I},\mathcal{J},\mathcal{K}\subset \Pi$, let \( \Gamma(\mathcal{I}, \mathcal{J}, \mathcal{K}) := \Phi_{\mathcal{I}} \cup (\Phi^+_\mathcal{J} \backslash \Phi_\mathcal{K}) \) be a subset of $\Phi$. In Proposition \ref{general closed condition}, we give a criteria in terms of $\mathcal{I},\mathcal{J}, \mathcal{K}$ on when $\Gamma(\mathcal{I}, \mathcal{J}, \mathcal{K})$ is closed, so that it corresponds to a subgroup \(H_{\Gamma(\mathcal{I}, \mathcal{J}, \mathcal{K})}\). The fixed point set \( \mathbb{X}_{\mathcal{L}, \Gamma(\mathcal{I}, \mathcal{J}, \mathcal{K})} \) of the torus action constructed above (Chapter \ref{Section Some fixed points subsets}) on \( T^*G/\!\!/H_{\Gamma(\mathcal{I}, \mathcal{J}, \mathcal{K})} \underset{[\mathfrak{g}^*/G]}{\times} T^*G/\!\!/_{{e}_\mathcal{L}} {U}_{e_\mathcal{L}} \) has the following combinatorial description related to the Weyl group of \( G \).

\begin{theorem}[Theorem \ref{fixed point set theorem}]$$\mathbb{X}_{\mathcal{L}, \Gamma(\mathcal{I}, \mathcal{J}, \mathcal{K})} = \coprod_{\substack{\mathcal{L}_1, \mathcal{L}_2, \mathcal{L}_3 \subset \mathcal{L} \\ \mathcal{L} = \mathcal{L}_1 \coprod \mathcal{L}_2 \coprod \mathcal{L}_3}} A(\mathcal{I}, \mathcal{L}_1, \mathcal{I}) \cap \bar{A}(\mathcal{I}, \mathcal{L}_2, \mathcal{J}) \cap \bar{B}(\mathcal{I}, \mathcal{L}_3, \mathcal{K}, \mathcal{I})$$
where
$$
\begin{aligned}
    A(\mathcal{I}, \mathcal{L}_1, \mathcal{I}) &= (W_{\mathcal{L}_1} \backslash W/W_\mathcal{I})^{\mathrm{free}}\\
    \bar{A}(\mathcal{I}, \mathcal{L}_2, \mathcal{J}) &= {W/W_\mathcal{I}} \times_{W/W_\mathcal{J}} (\overline{W_{\mathcal{L}_2}} \backslash W/W_\mathcal{J})^{\mathrm{free}}\\
    \bar{B}(\mathcal{I}, \mathcal{L}_3, \mathcal{K}, \mathcal{I}) &= (\overline{W_{\mathcal{L}_3}} \backslash W/W_\mathcal{I})^{\mathrm{free}} \cap (W/W_\mathcal{I} \times_{W/W_{\mathcal{I} \cup \mathcal{K}}}((W/W_{\mathcal{Y}})^{W_{\mathcal{L}_3}}/W_\mathcal{X})
\end{aligned}$$
and $\mathcal{X},\mathcal{Y}$ satisfying $\mathcal{X}\coprod \mathcal{Y}=\mathcal{I}\cup\mathcal{K}$ are subsets defined in Proposition \ref{general closed condition}.
\end{theorem}

For example, when $\mathcal{J}=\Pi, \mathcal{K}\subset \mathcal{I}$, we have 
$H_{\Gamma(\mathcal{I},\Pi,\mathcal{K})}=P_{\mathcal{I}}$. In this case $$\bar{A}(\mathcal{I},\mathcal{L}_2,\Pi)=
\begin{cases}
    \varnothing, \text{ if }\mathcal{L}_2\neq \varnothing;\quad\\
    W/W_{\mathcal{I}}, \text{ if }\mathcal{L}_2= \varnothing.\quad
\end{cases} \text{ and  }\quad
 \bar{B}(\mathcal{I},\mathcal{L}_3,\mathcal{K},\mathcal{I})=
\begin{cases}
    \varnothing, \text{ if }\mathcal{L}_3\neq \varnothing;\\
    W/W_{\mathcal{I}}, \text{ if }\mathcal{L}_3= \varnothing.
\end{cases}
$$
The fixed point set of \( T^*G/\!\!/P_\mathcal{I}\underset{[\mathfrak{g}^*/G]}{\times} T^*G/\!\!/_{{e}_\mathcal{L}} {U}_{e_\mathcal{L}} \) is thus the set $(W_{\mathcal{L}}\backslash W/W_\mathcal{I})^{\mathrm{free}}$,  which reproduced the fixed point set computed by Hoang, Krylov,  and Matvieievskyi in \cite{hoang2024around}.

Since the 4d branes in the second column of the list (\ref{list}) are all twisted cotangent bundles of \( G \) over a unipotent subgroup, one might naturally expect that the 4d mirror brane to 
$
T^*G / \!\!/ H_{\Gamma}
$
is always of the form 
$
T^*\langlands{G} / \!\!/_{\langlands{\psi}_{\Gamma}} \langlands{U}(\Gamma),
$ (which requires to take affinization when necessary,) where \( \langlands{U}(\Gamma) \) is certain unipotent subgroup of \( \langlands{G} \) that depends on \( \Gamma \) with a character \( \langlands{\psi}_{\Gamma} \). We expect that the fixed point set of 
\[
T^*G/\!\!/H_{\Gamma} \underset{[\mathfrak{g}^*/G]}{\times} T^*G/\!\!/_{{e}_\mathcal{L}} {U}_{e_\mathcal{L}}
\]
calculated in Lemma \ref{general fixed point} and Theorem \ref{fixed point set theorem} can be identified with the fixed point set of the intersection of 4d mirror branes, i.e., 
\[
T^*\langlands{G} / \!\!/_{\langlands{\psi}_{\Gamma}} \langlands{U}(\Gamma) \underset{[{\langlands{\mathfrak{g}}}^*/\langlands{G}]}{\times} T^*\langlands{G}/\!\!/ \langlands{P}_{\mathcal{L}}.
\]
\subsection*{Acknowledgments}
The authors thank David Ben-Zvi, Kifung Chan, Tom Gannon, Sanath Devalapurkar, Dmytro Matvieievsky, Ben Webster, Junzhe Lyu, Kaitao Xie, and Shilin Yu for valuable discussions
on various stages of this project. The work of N. C. Leung described in this paper was
substantially supported by grants from the Research Grants Council of the Hong Kong
Special Administrative Region, China (Project No. CUHK14306322, CUHK14305923 and CUHK14302224) and a direct grant from CUHK.
\section{Fundamental Conventions}
\subsection{Some Lie theoretic data}
\subsubsection{Algebraic groups and actions}
We assume every thing is over the complex number field $\mathbb{C}$. Let \( G \) be an algebraic group with Lie algebra \( \mathfrak{g} \). For a subgroup \( H \) of \( G \), we define the centralizer group of \( H \) in \( G \) as 
\[
C_G(H) := \{ g \in G \mid gh = hg, \forall h \in H \}.
\]

The normalizer group of \( H \) in \( G \) is given by 
\[
N_G(H) := \{ g \in G \mid gHg^{-1} = H \}.
\]

We denote by \( G^\circ \) the identity component of \( G \).

The adjoint and coadjoint representations of \( G \) are expressed as follows:
\[
\mathrm{Ad}: G \longrightarrow \operatorname{GL}(\mathfrak{g}) \quad \text{and} \quad \mathrm{Ad}^*: G \longrightarrow \mathrm{GL}\left(\mathfrak{g}^*\right).
\]

Similarly, the adjoint and coadjoint representations of \( \mathfrak{g} \) are defined by
\[
\mathrm{ad}: \mathfrak{g} \longrightarrow \mathfrak{gl}(\mathfrak{g}) \quad \text{and} \quad \mathrm{ad}^*: \mathfrak{g} \longrightarrow \mathfrak{gl}\left(\mathfrak{g}^*\right).
\]

For each \( y \in \mathfrak{g} \), we define the centralizer subalgebra as 
\[
C_\mathfrak{g}(y) := \{ x \in \mathfrak{g} \mid [x, y] = 0 \}.
\]

This subalgebra corresponds to the Lie algebra of the centralizer group 
\[
C_G(y) := \{ g \in G \mid \operatorname{Ad}_g(y) = y \}.
\]

When \( G \) acts on a variety \( X \), we denote the stabilizer subgroup of a point \( x \in X \) as 
\[
G_x := \{ g \in G \mid g \cdot x = x \}.
\]

We also define 
\[
X^G := \{ x \in X \mid G_x = G \}
\]

and 
\[
_GX := \{ x \in X \mid G_x = \{\operatorname{id}\} \},
\]

where \( \operatorname{id} \) is the identity element in \( G \).

\subsubsection{Semisimple groups}In the rest, we assume \( G \) is a connected semisimple algebraic group. We fix a choice of Borel subgroup \( B_G \) and a maximal torus \( T_G \subset B_G \). We will omit the subscript \( G \) when the context is clear. Let \( \Phi \) be the root system corresponding to the pair \( (G, B) \). Thus, \( \mathfrak{g} \) has a root space decomposition 
\[
\mathfrak{g} = \mathfrak{t} \oplus \bigoplus_{\alpha \in \Phi} \mathfrak{g}_\alpha
\]
with a choice of simple roots \( \alpha_1, \alpha_2, \ldots, \alpha_l \). Denote by \( \Pi \) the set of simple roots. We fix a root vector \( e_i \neq 0 \in \mathfrak{g}_{\alpha_i} \) and \( f_i \neq 0 \in \mathfrak{g}_{-\alpha_i} \). Let \( \Phi^+ \) (resp. \( \Phi^- \)) be the set of positive (resp. negative) roots. For \( \alpha \in \Phi \), we denote \( s_\alpha \) as the reflection associated with \( \alpha \). For \( \mathcal{I} \subset \Pi \), we will denote \( P_\mathcal{I} \) as the standard parabolic subgroup associated with \( \mathcal{I} \) containing \( B \), and \( L_\mathcal{I} \) be as the Levi subgroup of \( P_\mathcal{I} \). If \( P \) is a parabolic subgroup, then \( G/P \) is usually called the generalized flag variety.

Let \( W := N_{G}(T)/T \) be the Weyl group of \( G \). Then \( W \) is isomorphic to the group generated by the simple reflections \( s_{\alpha_1}, s_{\alpha_2}, \ldots, s_{\alpha_l} \).

Let \( A \subset \Phi \) be a subset of \( \Phi \). Let \( \Phi_A \) be the minimal root system containing \( A \). We will write \( -A := \{ \alpha \in \Phi \mid -\alpha \in A \} \). Then \( \Phi \backslash (-A) = -(\Phi \backslash A) \). Let \( A_s = \{ \alpha \in A \mid -\alpha \in A \} \) and \( A_u = A \backslash A_s \). Thus, \( A_s = A \cap (-A) \). Let \( V_A := \bigoplus_{\alpha \in A} \mathfrak{g}_\alpha \) be a subspace of \( \mathfrak{g} \) associated with \( A \).

Let \( W_A \) be the subgroup of \( W \) generated by the reflections associated with root vectors in \( A \). If \( I \) is a subset of \( \Pi \), then \( W_I \) is called a \textit{standard parabolic subgroup} of \( W \), and conjugate subgroups \( wW_Iw^{-1}, w \in W \) are called \textit{parabolic subgroups} of \( W \).

Let \( \langle \cdot, \cdot \rangle: \mathfrak{g} \otimes_\mathbb{C} \mathfrak{g} \longrightarrow \mathbb{C} \) be the Killing form \( \kappa \) defined by
\[
\langle \xi, \eta \rangle = \operatorname{trace}\left(\operatorname{ad}_{\xi} \circ \operatorname{ad}_\eta\right)
\]
for all \( \xi, \eta \in \mathfrak{g} \). The Killing form is \( G \)-invariant, symmetric, and non-degenerate. Let \( (\cdot, \cdot) \) be the inner product of the root system \( \Phi \).

Given a vector subspace \( V \subset \mathfrak{g} \), we write
\[
V^{\perp} := \{ \xi \in \mathfrak{g} : \langle \xi, \eta \rangle = 0 \text{ for all } \eta \in V \}.
\]

Then \( (V^\perp)^\perp = V \). It is easy to see that \( \mathfrak{g}_\alpha^\perp = \mathfrak{t} \oplus \bigoplus_{\beta \in \Phi \backslash \{-\alpha\}} \mathfrak{g}_\beta \) and \( \mathfrak{t}^\perp = \bigoplus_{\alpha \in \Phi} \mathfrak{g}_\alpha \). There is a natural isomorphism \( (\mathfrak{g}/V)^* \cong V^\perp \) via the Killing form.

Let \( \mathcal{N} \) be the nilpotent cone of \( \mathfrak{g} \). Let \( U \) be the unipotent radical of \( B \).

Let \( e \) be a nilpotent element of \( \mathfrak{g} \). Recall that \( \tau = (e, h, f) \in \mathfrak{g}^{\oplus 3} \) is an \( \mathfrak{sl}_2 \)-triple if 
\[
[e, f] = h, \quad [h, e] = 2e, \quad [h, f] = -2f.
\]

The associated \textit{Slodowy slice} is defined by 
\[
S_\tau := e + C_\mathfrak{g}(f) \subset \mathfrak{g}.
\]

Fix \( e \in \mathcal{N} \); then any two \( \mathfrak{sl}_2 \)-triples \( (e, h_1, f_1), (e, h_2, f_2) \) are conjugate by an element in \( C_G(e) \). Therefore, \( S_\tau \) is uniquely determined by the nilpotent orbit of \( e \) up to conjugation, and we shall write \( S_e \) for \( S_\tau \) when the context is clear. Then \( S_{-e} \) is the Slodowy slice associated with the \( \mathfrak{sl}_2 \)-triple \( (-e, h, -f) \), and \( S_f \) is the Slodowy slice associated with the \( \mathfrak{sl}_2 \)-triple \( (-f, -h, e) \).

\subsubsection{Hamiltonian action on the cotangent bundle of $G$}
We will be interested in the case \( X = T^*G = G \times \mathfrak{g}^* \), which is a symplectic Hamiltonian \( G \times G \)-variety obtained by lifting the left and right multiplications on \( G \).

Through the Killing form, we can identify \( \mathfrak{g} \) and \( \mathfrak{g}^* \) via the map
\[
\begin{aligned}
\kappa: \mathfrak{g} & \longrightarrow \mathfrak{g}^* \\
x & \longmapsto \langle x, - \rangle.
\end{aligned}
\]

The action of \( G \times G \) on \( T^*G \cong G \times \mathfrak{g} \) is given by
\[
(g, h) \cdot (a, x) = (gah^{-1}, \operatorname{Ad}_h x), \quad \text{for } g, h \in G, (a, x) \in G \times \mathfrak{g}.
\]

The associated moment map is then 
\[
\begin{aligned}
\mu: T^* G \cong G \times \mathfrak{g} & \longrightarrow \mathfrak{g} \times \mathfrak{g} \\
(a, x) & \longmapsto \left( \operatorname{Ad}_a x, x \right).
\end{aligned}
\]

Let \( \mu_L, \mu_R \) be the moment maps that arise from \( \mu \) by projecting to the first and second factors, respectively. Then \( \mu_L \) (resp. \( \mu_R \)) corresponds to the moment map of the left (resp. right) action of \( G \). 

There is an involution \( \Theta: G \times \mathfrak{g} \rightarrow G \times \mathfrak{g} \) defined by 
\[
\Theta(a, x) = (a^{-1}, \operatorname{Ad}_{a} x),
\]
which is compatible with \( \mu \) and switches the left and right actions. More precisely, let \( \theta: \mathfrak{g} \times \mathfrak{g} \rightarrow \mathfrak{g} \times \mathfrak{g} \) be defined by \( \theta(x, y) = (y, x) \). Then we have 
\[
\mu \circ \Theta = \theta \circ \mu
\]
and 
\[
\Theta((g, h) \cdot (a, x)) = (h, g) \cdot \Theta(a, x).
\]

\subsubsection{Algebraic symplectic quotients}Let \( G_\eta \subset G \) be the \( G \)-stabilizer of \( \eta \in \mathfrak{g}^* \) under the coadjoint action. Let \( X \) be a symplectic Hamiltonian \( G \)-variety with moment map \( \mu: X \rightarrow \mathfrak{g}^* \). Then \( \mu^{-1}(\eta) \subset X \) is \( G_\eta \)-invariant.

We define 
\[
X/\!\!/_\eta G := \mu^{-1}(\eta)/G_\eta.
\]

When \( G_\eta \) acts freely on \( \mu^{-1}(\eta) \) and this action admits a geometric quotient \( \mu^{-1}(\eta) \rightarrow X/\!\!/_\eta G \), the quotient variety \( X/\!\!/_\eta G \) is smooth and inherits a unique algebraic symplectic form from \( X \). When \( \eta = 0 \), we will write \( X/\!\!/_0 G \) simply as \( X/\!\!/G \).

\subsubsection{Associated bundle}
Let \( H \) be a closed subgroup of \( G \) and \( Y \) be an \( H \)-variety. The space \( G \times_H Y := (G \times Y)/H \) is a fiber bundle over \( G/H \) such that each fiber is isomorphic to \( Y \). It is called the associated bundle for \( G \), \( H \), and \( Y \). An element in \( G \times_H Y \) can be denoted as \( [g, y] \), where \( g \in G \) and \( y \in Y \).

The following lemma will be useful later.

\begin{lemma}\label{embedding lemma}
Let \( G \) be an algebraic group and \( H \) be an algebraic subgroup of \( G \). Let \( Y \) be an \( H \)-invariant subvariety of \( G \)-representation \( V \). Then 
\[
\begin{aligned}
    \iota: G \times_H Y & \rightarrow (G/H) \times V \\
    [g, v] & \mapsto (gH, gv)
\end{aligned}
\]
is a \( G \)-equivariant embedding.
\end{lemma}

\begin{prf}
Since $G\times Y\rightarrow (G/H)\times V, (g,v)\mapsto (gH,gv)$ is independent of the $H$-action, $\iota$ is indeed a morphism. It is straightforward to show that it is equivariant. We prove that it is an embedding as follows. If \( (g_1H, g_1 v_1) = (g_2H, g_2 v_2) \), then \( g_1 = g_2 h \) for some \( h \in H \). Furthermore, \( g_1 v_1 = g_2 v_2 \) implies \( g_2 v_1 = g_2 v_2 \) which leads to \( v_2 = h v_1 \). Therefore, \( [g_1, v_1] = [g_2, v_2] \in G \times_H V \).
\end{prf}

\subsection{Twisted Cotangent Bundles}\label{twisted cotangent bundles}
The notations and formulations in this subsection mostly follow \cite{crooks2023universal}.

\subsubsection{Affine Hamiltonian Lagrangian (AHL) \( G \)-bundles}

Let \( X \) be a symplectic variety with symplectic form \( \omega \). A fiber bundle \( (X, Y, \sigma) \), with \( \sigma: X \rightarrow Y \), is a \textit{Lagrangian fibration} if \( \sigma^{-1}(y) \) is a Lagrangian subvariety in \( X \) for any \( y \in Y \).

We denote the canonical vector bundle isomorphism induced by \( \omega \) as 
\[
T^* X \rightarrow T X: v \mapsto v^{\#}.
\]

The fibers of $\sigma$ carry an additional structure given by 
\[
T^* Y \times_Y X \rightarrow T X: (p, x) \mapsto \xi^p(x) = \left( \left( T_x \sigma \right)^*(p) \right)^{\#}.
\]

We call the Lagrangian fibration \( \sigma: X \rightarrow Y \) \textit{affine} if all fibers of \( \sigma \) are simply-connected and all vector fields \( \xi^p \) are complete.

If, moreover, \( X \) is a Hamiltonian \( G \)-variety with moment map \( \mu \), \( Y \) is a \( G \)-variety, and the affine Lagrangian fibration \( \sigma: X \rightarrow Y \) is \( G \)-equivariant, then \( (X, \mu, \sigma) \) is called an \textit{affine Hamiltonian Lagrangian (AHL) \( G \)-bundle}.

Two AHL \( G \)-bundles \( (X_1, \mu_1, \sigma_1) \) and \( (X_2, \mu_2, \sigma_2) \) over \( Y \) are called isomorphic if there exists a \( G \)-equivariant symplectic isomorphism \( \varphi: X_1 \longrightarrow X_2 \) satisfying \( \mu_1 = \mu_2 \circ \varphi \) and \( \sigma_1 = \sigma_2 \circ \varphi \). We denote by \( \mathscr{AHL}_{G}(Y) \) the set of isomorphism classes of AHL \( G \)-bundles over \( Y \).

\subsubsection{AHL \( G \)-bundles over homogeneous spaces}

Let \( H \subset G \) be a closed subgroup with Lie algebra \( \mathfrak{h} \subset \mathfrak{g} \). Let \( \psi \in (\mathfrak{h}^*)^H \).

The subgroup \( H \) of \( G \) acts on \( T^*G = G \times \mathfrak{g} \) from the left, viewing it as the subgroup \( H \times \{\operatorname{id}\} \subset G \times G \) (resp. from the right, viewing it as the subgroup \( \{\operatorname{id}\} \times H \subset G \times G \)), with moment maps \( \mu_{L,H} \) (resp. \( \mu_{R,H} \)): 
\[
G \times \mathfrak{g} \rightarrow \mathfrak{h}^*.
\]

We define 
\[
T^*G/\!\!/_\psi H := \mu_{R,H}^{-1}(\psi)/H
\]
and 
\[
H_\psi \backslash\!\!\backslash T^*G := H \backslash \mu_{L,H}^{-1}(\psi).
\]

Let \( \hat{\psi} \) be a lift of \( \psi \) in \( \mathfrak{g} \) under the identification by the Killing form. Then 
\[
\mu_{R,H}^{-1}(\psi) = \{ (a, x) \mid x \in \hat{\psi} + \mathfrak{h}^\perp \}
\]
and 
\[
\mu_{L,H}^{-1}(\psi) = \{ (a, x) \mid \operatorname{Ad}_a x \in \hat{\psi} + \mathfrak{h}^\perp \}.
\]

The involution \( \Theta \) on $T^*G$ restricts to an isomorphism between \( \mu_{R,H}^{-1}(\psi) \) and \( \mu_{L,H}^{-1}(\psi) \), and it also descends to an isomorphism between \( T^*G/\!\!/_\psi H \) and \( H_\psi \backslash\!\!\backslash T^*G \). The inverse is also given by \( \Theta \).

If \( g \in N_G(H) \), then \( \operatorname{Ad}_g \) preserves \( \mathfrak{h} \) and hence \( N_G(H) \) acts on \( \mathfrak{h}^* \). Let \( N_G(H)_{\psi} \) be the subgroup of \( N_G(H) \) that fixes \( \psi \). The space \( T^*G/\!\!/_\psi H \) admits a \( G \times N_G(H)_{\psi} \)-Hamiltonian action.

We have 
\[
T^*G/\!\!/_\psi H = G \times_H (\hat{\psi} + \mathfrak{h}^\perp)
\]
which is also the associated bundle for \( G, H \), and \( \hat{\psi} + \mathfrak{h}^\perp \). The moment map under the \( G \)-action is 
\[
G \times_H (\hat{\psi} + \mathfrak{h}^\perp) \rightarrow \mathfrak{g}, \quad [g, x] \mapsto \operatorname{Ad}_g x.
\]

The projection to \( G/H \) is an affine Lagrangian fibration. Thus, the space \( T^*G/\!\!/_\psi H \) is an AHL \( G \)-bundle over \( G/H \). It is called the \( \psi \)-\textit{twisted cotangent bundle} of \( G/H \). When \( \psi = 0 \), we have 
\[
T^*G/\!\!/ H = G \times_H \mathfrak{h}^\perp = T^*(G/H)
\]
which is exactly the cotangent bundle of \( G/H \).

In \cite{bisiecki1985holomorphic}, the author classifies all AHL \( G \)-bundles over \( G/H \) and finds that they are parametrized by \( (\mathfrak{h}^*)^H \).

\begin{proposition}[\cite{bisiecki1985holomorphic}]
The map
\[
\left( \mathfrak{h}^* \right)^H \longrightarrow \mathscr{AHL}_G(G/H), \quad \psi \mapsto \left[T^*G/\!\!/_\psi H\right]
\]
is a bijection.
\end{proposition}

\subsection{Symplectic duality and S-duality}
Fix a reductive group \( G \) (the “gauge group”) and a representation \( \mathrm{N} \) (the “matter”). This data defines in physics a 3d \( \mathcal{N} = 4 \) supersymmetric gauge theory \cite{bullimore2017coulomb,bullimore2016boundaries,nakajima2015towards}. It has various moduli of vacua, one of which is the Higgs branch, which is the Hamiltonian reduction of \( T^*\mathrm{N} \) by $G$. Another moduli of vacua is called the Coulomb branch, denoted \( M_C(G, \mathrm{N}) \), which is rigorously defined in \cite{braverman2016towards}. The Higgs and Coulomb branches are expected to be symplectic dual to each other. Many of these moduli of vacua are AHL $G$-bundles.

The current hypothesis states that \( \mathrm{N} \) is a smooth affine algebraic variety, as mentioned in \cite{Nakajima:2024mlb}. Define $\mathcal{F}=\mathbb{C}((z))$ and $\mathcal{O}=\mathbb{C}[[z]]$, representing the field of formal Laurent series and the ring of formal power series, respectively. Consider $\operatorname{Gr}_G:=G(\mathcal{F})/G(\mathcal{O})$. They then described the variety of triples $R_{G,\mathrm{N}}$ as:
$$R_{G,\mathrm{N}}=\{[g(z),s(z)]\in G(\mathcal{F})\times^{G(\mathcal{O})}\mathrm{N}(\mathcal{O})\mid g(z)s(z)\in \mathrm{N}(\mathcal{O})\}.$$
The group $G(\mathcal{O})$ acts on $R_{G,\mathrm{N}}$ from the left. The Coulomb branch is characterized by:
$$M_C(G,\mathrm{N}):=\operatorname{Spec}H_*^{G(\mathcal{O})}(R_{G,\mathrm{N}}).$$

Define \( D_{G(\mathcal{O})}(\operatorname{Gr}_G) \) as the derived category of \( G(\mathcal{O}) \)-equivariant sheaves on \( \operatorname{Gr}_G \), which is equipped with a natural monoidal convolution structure denoted by \( * \). Let \( \omega_{R_{G,\mathrm{N}}} \) represent the dualizing sheaf of \( R_{G,\mathrm{N}} \). Denote by 
\(
\mathcal{A}_\mathrm{N} := \pi_*(\omega_{R_{G,\mathrm{N}}})
\)
the pushforward object in \( D_{G(\mathcal{O})}(\operatorname{Gr}_G) \). Then, we naturally have the equality:
\[
H_*^{G(\mathcal{O})}(R_{G,\mathrm{N}}) = H_*^{G(\mathcal{O})}(\operatorname{Gr}_G, \mathcal{A}_\mathrm{N}).
\]

$\mathcal{A}_\mathrm{N}$ is a commutative ring object in $D_{G(\mathcal{O})}(\operatorname{Gr}_G)$ equipped with the homomorphism $\mathcal{A}_\mathrm{N}\times \mathcal{A}_\mathrm{N}\rightarrow \mathcal{A}_\mathrm{N}$.

A provisional definition of the S-dual for the action \( G \curvearrowright T^* \mathrm{N} \) can be found in \cite{Nakajima:2024mlb}. Here we briefly revisited it.

Recall that the geometric Satake equivalence establishes an equivalence of tensor categories:
\[
\left( \operatorname{Perv}_{G(\mathcal{O})} \left( \operatorname{Gr}_G \right), * \right) \cong \left( \operatorname{Rep} \langlands{G}, \otimes \right),
\]
where the left-hand side is the category of \( G(\mathcal{O}) \)-equivariant perverse sheaves on \( \operatorname{Gr}_G \), and the right-hand side is the tensor category of finite-dimensional representations of the Langlands dual group \( \langlands{G} \).

Let \( \mathbb{C}[\langlands{G}] \) denote the regular representation of \( \langlands{G} \) appearing on the right-hand side. Its corresponding object on the left-hand side is referred to as the \textit{regular sheaf} and is denoted by \( \mathcal{A}_{\mathbb{C}[\langlands{G}]} \).

The S-dual associated with \( G \curvearrowright T^* \mathrm{N} \) is given by:
\[
\langlands{\mathrm{M}} := \operatorname{Spec} H_*^{G(\mathcal{O})} \left( \operatorname{Gr}_G, \mathcal{A}_\mathrm{N} \otimes^! \mathcal{A}_{\mathbb{C}[\langlands{G}]} \right).
\]

It was proposed that \( T^*G /\!\!/ G \) corresponds to the mirror 4d brane of 
\(
T^*\langlands{G} /\!\!/_{\langlands{e_\Pi}} \langlands{U} = \langlands{G} \times S_{\langlands{e_\Pi}}.
\)
Consequently, \( T^* \mathrm{N} /\!\!/ G \) is expected to be 3d mirror dual to 
\(
\langlands{\mathrm{M}} /\!\!/_{\langlands{e_\Pi}} \langlands{U}.
\)
A rough proof of this statement, as communicated to the author, was provided by Sanath Devalapurkar.
\begin{proposition}
\label{proposition Sympl dual}
    The coordinate ring of $\langlands{\mathrm{M}}/\!\!/_{\langlands{e_\Pi}} \langlands{U}$ is $H_*^{G(\mathcal{O})}(\operatorname{Gr}_G, \mathcal{A}_\mathrm{N})$.
\end{proposition}
\begin{prf}
    We have the derived Satake equivalence 
    $$D_{G(\mathcal{O})}(\operatorname{Gr}_G)\simeq \operatorname{QCoh}([(\langlands{\mathfrak{g}})^*/\langlands{G}])$$

    Then the image of the ring object $\mathcal{A}_\mathrm{N}$ in $\operatorname{QCoh}([(\langlands{\mathfrak{g}})^*/\langlands{G}])$ is the pushforward of the structure sheaf along the 1-shifted Lagrangian $\pi: [\langlands{M}/\langlands{G}]\rightarrow [(\langlands{\mathfrak{g}})^*/\langlands{G}]$, say $\pi_*\mathcal{O}_{[\langlands{M}/\langlands{G}]}$.

    The functor $p_*:D_{G(\mathcal{O})}(\operatorname{Gr}_G) \rightarrow D_{G(\mathcal{O})}(pt)$ given by cohomology pushforward to a point identifies under derived  Satake with the functor $\kappa^*: \operatorname{QCoh}([(\langlands{\mathfrak{g}})^*/\langlands{G}])\rightarrow \operatorname{QCoh}((\langlands{\mathfrak{g}})^*/\!\!/\langlands{G})$ given by pullback along the Kostant slice $(\langlands{\mathfrak{g}})^*/\!\!/\langlands{G}\rightarrow (\langlands{\mathfrak{g}})^*/\langlands{G}$. Taking cohomology of both sides, we have 
    $$H_*^{G(\mathcal{O})}(\operatorname{Gr}_G, \mathcal{A}_\mathrm{N})\cong \Gamma(\langlands{\mathfrak{g}})^*/\!\!/\langlands{G}, \kappa^*\pi_*\mathcal{O}_{[\langlands{M}/\langlands{G}]})\cong \mathbb{C}[\langlands{M}\times_{(\langlands{\mathfrak{g}})^*}S_{\langlands{e_\Pi}}].$$    
\end{prf}

\section{Intersection of Hamiltonian $G$-varieties}\label{Intersection of Hamiltonian $G$-varieties}
Suppose $X_i,i=1,2$ are $G$-Hamiltonian varieties with moment map $\mu_i:X_i\rightarrow\mathfrak{g}^*$. We denote  $\underline{X}_i=(G\curvearrowright X\stackrel{\mu_i}{\longrightarrow}\mathfrak{g}^*)$. The intersection of \( \underline{X}_1 \) and \( \underline{X}_2 \) is defined as the holomorphic symplectic quotient 
\[
\begin{aligned}
    \underline{X}_1 \underset{[\mathfrak{g}^*/G]}{\times} \underline{X}_2 &:= \{(x_1,x_2) \in X_1  \times X_2 \mid \mu_{1}(x_1) = \mu_{2}(x_2) \}/G\\&=(X_1\times_{\mathfrak{g}^*}X_2)/G\\
    &=(X_1^-\times X_2)/\!\!/G.
\end{aligned}
\]

We are mostly interested in the intersections of twisted cotangent bundles. Many of them are symplectic resolutions of their affinizations, as mentioned in the introduction. 

A quick exercise is the following:
\begin{proposition}\label{turn intersection to double quotient}
For any subgroups $H_1,H_2$ in $G$ and characters $\psi_1\in (\mathfrak{h}_1)^{H_1}$ and $\psi_2\in (\mathfrak{h}_2)^{H_2}$, we have
    \( H_1  {}_{\psi_1}\!\! \backslash\!\!\backslash T^*G \underset{[\mathfrak{g}^*/G]}{\times} T^*G/\!\!/_{\psi_2} H_2 \cong T^*G/\!\!/_{(\psi_1, \psi_2)}(H_1 \times H_2)\)
\end{proposition}

\begin{prf}
We have
\[
\begin{aligned}
    H_1 {}_{\psi_1} \backslash\!\!\backslash T^*G &= H_1 \backslash \{(g_1, x_1) \in G \times \mathfrak{g} \mid \operatorname{Ad}_{g_1} x_1 \in \hat{\psi}_1 + \mathfrak{h}_1^\perp\} \\
    T^*G/\!\!/_{\psi_2} H_2 &= \{(g_2, x_2) \in G \times \mathfrak{g} \mid x_2 \in \hat{\psi}_2 + \mathfrak{h}_2^\perp\}/H_2
\end{aligned}
\]
Then, taking intersection means we impose the new condition \( x_1 = \operatorname{Ad}_{g_2} x_2  \) and then take the quotient by \( G \). In \( H_1 {}_{\psi_1} \backslash\!\!\backslash T^*G \underset{[\mathfrak{g}^*/G]}{\times} T^*G/\!\!/_{\psi_2} H_2 \), the element \( [(g_1, x_1), (g_2, x_2)] \) is equivalent to \( [(g_1 g_2, x_2), (1, x_2)] \), and the isomorphism sends it to \( (g_1 g_2, x_2) \in T^*G/\!\!/_{(\psi_1, \psi_2)}(H_1 \times H_2) \).
\end{prf}

Via the involution \( \Theta \), we will write \( H_1 {}_{\psi_1}\!\! \backslash\!\!\backslash T^*G \underset{[\mathfrak{g}^*/G]}{\times} T^*G/\!\!/_{\psi_2} H_2 \) instead of \( T^*G/\!\!/_{\psi_1} H_1 \underset{[\mathfrak{g}^*/G]}{\times} T^*G/\!\!/_{\psi_2} H_2 \) for the beauty of symmetry.

\subsection{Intersection with Bielawski hyperkähler slices}\label{section bielawski slices}
Let \( e \) be a nilpotent element in \( \mathfrak{g} \) and $S_e$ be its associated Slodowy slice. Then \( G \times S_e \) is a hyperkähler and Hamiltonian \( G \)-variety studied by Bielawski \cite{bielawski1997hyperkahler, bielawski2017slices}, Moore–Tachikawa \cite{moore20122d}, and others. We recall its basic properties here. Suppose \( \tau = (e, h, f) \) is an \( \mathfrak{sl}_2 \)-triple associated with \( e \) and let \( \mathfrak{g} = \bigoplus \mathfrak{g}_k \) be the decomposition of \( \mathfrak{g} \) into eigenspaces for \( h \).

\begin{proposition}
    The space \( \mathfrak{g}_{-1} \) has a non-degenerate skew-symmetric bilinear form defined by \( \langle[-,-], e\rangle \).
\end{proposition}

\begin{prf}
Suppose \( \langle[x,y], e\rangle = 0 \) for any \( y \in \mathfrak{g}_{-1} \).
By associativity and non-degeneracy of the Killing form, we have \( [x, e] = 0 \).
Now \( C_\mathfrak{g}(e) \subset \mathfrak{g}_{\geq 0} := \bigoplus_{k \geq 0} \mathfrak{g}_k \) implies \( x = 0 \).
\end{prf}

Let \( \mathfrak{l} \subset \mathfrak{g}_{-1} \) be a Lagrangian subspace with respect to this form.  
Let \( \mathfrak{u}_e := \mathfrak{l} \oplus \bigoplus_{k \leq -2} \mathfrak{g}_k \).
Let \( U_e \subset G \) be the associated unipotent group.

\begin{example}Suppose \(\mathfrak{g}\) is of type \(A_4\). Let \(e = E_{12} + E_{34} + E_{45}\). Then,
\[
\mathfrak{g}_{-1} = \left\{ 
\begin{pmatrix}
0 & 0 & x & 0 & 0 \\ 
0 & 0 & 0 & y & 0 \\ 
0 & 0 & 0 & 0 & 0 \\ 
z & 0 & 0 & 0 & 0 \\ 
0 & w & 0 & 0 & 0 
\end{pmatrix} \mid x, y, z, w \in \mathbb{C}
\right\}
\]
A possible choice for a Lagrangian subspace of \(\mathfrak{g}_{-1}\) consists of the matrices with \(x = y = 0\). Thus, elements of \(\mathfrak{u}_{e} = \mathfrak{l} \oplus \bigoplus_{k \leq -2} \mathfrak{g}_k\) can take the following form:
\[
\begin{pmatrix}
0 & 0 & 0 & 0 & 0 \\ 
* & 0 & * & 0 & 0 \\ 
0 & 0 & 0 & 0 & 0 \\ 
* & 0 & * & 0 & 0 \\ 
* & * & * & * & 0 
\end{pmatrix}
\]
\end{example}

We have a natural character \( \chi_e := \langle e, - \rangle: \mathfrak{u}_e \rightarrow \mathbb{C} \).
\begin{lemma}
    \( \chi_e \in (\mathfrak{u}_e^*)^{U_e} \).
\end{lemma}

\begin{prf}
    It suffices to show that \( \langle e, [\mathfrak{u}_e, \mathfrak{u}_e] \rangle = 0 \). Since \( \mathfrak{l} \) is a Lagrangian subspace, it suffices to show \( \langle e, [\mathfrak{u}_e, \mathfrak{g}_{\leq -2}] \rangle = \langle \mathfrak{u}_e, [e, \mathfrak{g}_{\leq -2}] \rangle = 0 \). This is because \( \langle \mathfrak{g}_{\leq -1}, \mathfrak{g}_{\leq 0} \rangle = 0 \).
\end{prf}

Let \( \pi_e: \mathfrak{g} \rightarrow \mathfrak{u}_e^* \) be induced from \( \mathfrak{g}^* \rightarrow \mathfrak{u}_e^* \) and the Killing form identification. Then 
\[
\pi_e^{-1}(\chi_e) = e + \mathfrak{u}_e^{\perp}.
\]

\begin{proposition}[\cite{gan2002quantization}]\label{adjoint action map}
    The adjoint action map \( \alpha_e: U_e \times S_e \rightarrow e + \mathfrak{u}_e^{\perp}, (u,x) \mapsto \operatorname{Ad}_u x \) is an isomorphism of affine varieties.
\end{proposition}

As a result, we have \( G \times S_e \cong T^*G/\!\!/_{\chi_e} U_e \). We shall write \( T^*G/\!\!/_{e} U_e \) for convenience when the context is clear.

\begin{proposition}
    Let \( X \) be a Hamiltonian \( G \)-variety with moment map \( \varrho: X \rightarrow \mathfrak{g} \).  The intersection \( \underline{X} \underset{[\mathfrak{g}^*/G]}{\times} T^*G/\!\!/_e U_e \) is isomorphic to $
\varrho^{-1}(S_{e}) = X \times_\mathfrak{g} S_{e}.$
\end{proposition}

\begin{prf}
    If \( (x, g, y) \in X \times G \times S_e \), we can apply the action of \( G \) to assume \( g = \operatorname{id} \). Then the condition \( \varrho(x) = y \) implies \( x \in \varrho^{-1}(S_{e}) \). The inverse map from \( \varrho^{-1}(S_{e}) \) is the inclusion \( x \mapsto (x, 1, \varrho(x)) \in X \underset{[\mathfrak{g}^*/G]}{\times} T^*G/\!\!/_e U_e \).
\end{prf}

\subsection{Universal centralizer and intersections of AHL $G$-bundles over $G/U$}The elements \( e_{\Pi} := \sum_{i \in \Pi} e_i \) and \( f_{\Pi} = \sum_{i \in \Pi} f_i \) are regular nilpotent elements (see for example, \cite[Corollary 3.2.13]{chriss1997representation}) and \( \mathcal{S} := S_{e_{\Pi}} \cong S_{f_\Pi} \) is called the Kostant slice. It is a regular section of the flat morphism \( \mathfrak{g} \rightarrow \mathfrak{g}/\!\!/G \), where \( \mathfrak{g}/\!\!/G := \operatorname{Spec} \mathbb{C}[\mathfrak{g}]^G \) parametrizes the set of closed adjoint orbits.

For $\mathcal{I}\subset \Pi$, let \( f_{\mathcal{I}} := \sum_{i \in \mathcal{I}} f_i \). Then \( \chi_{f_{\mathcal{I}}} := \langle f_\mathcal{I}, - \rangle: \mathfrak{u} \rightarrow \mathbb{C} \) is a character such that \( \chi_{f_{\mathcal{I}}} \in (\mathfrak{u}^*)^U \). In fact, \( \langle f_{\mathcal{I}}, [\mathfrak{u}, \mathfrak{u}] \rangle = \langle [f_{\mathcal{I}}, \mathfrak{u}], \mathfrak{u} \rangle = 0 \) as \( [f_{\mathcal{I}}, \mathfrak{u}] \subset \mathfrak{b} \). Let \( T^*G/\!\!/_{{f_{\mathcal{I}}}} U := T^*G/\!\!/_{\chi_{f_{\mathcal{I}}}} U \). Then 
\[
\begin{aligned}
Y_{\mathcal{I}} &:= U {}_{f_{\mathcal{I}}}\!\! \backslash\!\!\backslash T^*G \underset{[\mathfrak{g}^*/G]}{\times} T^*G/\!\!/_{f_{\Pi}} U \\
&= U \backslash \{(g,x) \in G \times \mathfrak{g} \mid x \in f_\Pi + \mathfrak{b}, \operatorname{Ad}_g x \in f_{\mathcal{I}} + \mathfrak{b}\}/U \\
&= \{([g],x) \in U \backslash G \times \mathcal{S} \mid \operatorname{Ad}_g x \in f_{\mathcal{I}} + \mathfrak{b}\},
\end{aligned}
\]
where the second equation is obtained by Proposition \ref{adjoint action map} for the special case when \( e \) is regular, and \( U \) is the unipotent radical of a Borel subgroup.

The \textit{universal centralizer} of \( G \) is the affine variety 
\[
J_G := \{(g,x) \in G \times \mathfrak{g} \mid \operatorname{Ad}_g(x) = x\}.
\]

It plays an important role in geometric representation theory and the theory of symplectic duality. One can also interpret it as the intersection of Hamiltonian \( G \)-varieties. In fact, we have \cite{balibanu2017partial,jin2022homological}

\begin{proposition}\label{UniversalCent as intersection}
    \( J_G \) is isomorphic to \( U {}_{f_{\Pi}} \!\!\backslash\!\!\backslash T^*G \underset{[\mathfrak{g}^*/G]}{\times} T^*G/\!\!/_{f_{\Pi}} U \).
\end{proposition}

\begin{prf}
    According to the above discussion, it suffices to show \( J_G \cong Y_{\Pi} = \{([g],x) \in U \backslash G \times \mathcal{S} \mid \operatorname{Ad}_g x \in f_{\Pi} + \mathfrak{b}\} \). By Proposition \ref{adjoint action map}, there is a unique element \( g_0 \) in the coset \( Ug \) such that \( \operatorname{Ad}_{g_0} x \in \mathcal{S} \). Then \( g_0 \in C_G(x) \), and thus we define the isomorphism between \( J_G \) and \( Y_\Pi \).
\end{prf}

Let \( \langlands{G} \) be the Langlands dual group of \( G \). Recall that \( J_{\langlands{G}} \) is the Coulomb branch of a 3d \( \mathcal{N} = 4 \) SUSY gauge theory associated with \( (G,0) \), where \( '0' \) is the trivial 0-dimensional representation of \( G \). The Higgs branch of \( (G,0) \) could also be written as an intersection between twisted cotangent bundles \( G \backslash\!\!\backslash T^*G \underset{[\mathfrak{g}^*/G]}{\times} T^*G/\!\!/G \). We will also see how 3d mirror dual bow varieties interpret this in Section \ref{section bow variety}.

\subsection{Intersection of $P\backslash\!\!\backslash T^*G$ and AHL $G$-bundle over $G/U$}Let \( \mathcal{I}, \mathcal{L} \subset \Pi \).
Suppose \( P_{\mathcal{I}} \) is a standard parabolic subgroup associated to \( \mathcal{I} \), and \( f_\mathcal{L} := \sum_{i \in \mathcal{L}} f_i \). Let 
\[
Z_{\mathcal{I}, \mathcal{L}} := P_\mathcal{I} \backslash\!\!\backslash T^*G \underset{[\mathfrak{g}^*/G]}{\times} T^*G/\!\!/_{f_{\mathcal{L}}} U.
\]

Then 
\[
Z_{\mathcal{I}, \mathcal{L}} := P_\mathcal{I} \backslash \{(g,x) \in G \times \mathfrak{g} \mid \operatorname{Ad}_g x \in \mathfrak{p}_{\mathcal{I}}^\perp, x \in f_\mathcal{L} + \mathfrak{b}\}/U.
\]

It has a projection map \( \pi \) to the double coset space \( P_\mathcal{I} \backslash G/U \cong W_\mathcal{I} \backslash W \), $\pi:Z_{\mathcal{I}, \mathcal{L}}\rightarrow W_\mathcal{I} \backslash W$. For \( w \in W_\mathcal{I} \backslash W \), there exists \( x \in f_\mathcal{L} + \mathfrak{b} \) such that \( \operatorname{Ad}_w x \in \mathfrak{p}_{\mathcal{I}}^\perp \) if and only if \( \operatorname{Ad}_w f_\mathcal{L} \in \mathfrak{p}_{\mathcal{I}}^\perp \), and if and only if \( w(-\alpha_i) \in \Phi^+ \backslash \Phi_{\mathcal{I}} \) for \( i \in \mathcal{L} \). Then the image of the projection map is 
\[
\{w \in W_\mathcal{I} \backslash W \mid w(-\alpha_i) \in \Phi^+ \backslash \Phi_{\mathcal{I}}, i \in \mathcal{L}\}.
\]

In subsection \ref{Parabolic Slodowy variety}, we will see that this set is naturally bijective to \( (W_\mathcal{L} \backslash W/W_\mathcal{I})^{\operatorname{free}} \), the set of free \( W_\mathcal{L} \times W_\mathcal{I} \)-orbits in \( W \).

Suppose \( w \in G \) is a representative of its image in \( (W_\mathcal{L} \backslash W/W_\mathcal{I})^{\operatorname{free}} \). Then for \( x_1, x_2 \in f_\mathcal{L} + \mathfrak{b} \cap \operatorname{Ad}_{w^{-1}}(\mathfrak{p}_{\mathcal{I}}^\perp) \), we have \( (w, x_1) = (w, x_2) \in Z_{\mathcal{I}, \mathcal{L}} \) if and only if there exists \( p \in P_\mathcal{I}, u \in U \) such that \( pwu = w \) and \( \operatorname{Ad}_{u^{-1}} x_1 = x_2 \). Therefore, \( u = w^{-1} p w \in U \cap w^{-1} P_\mathcal{I} w \), and the fiber is the set of \( U \cap w^{-1} P_\mathcal{I} w \)-orbits in \( f_\mathcal{L} + \mathfrak{b} \cap \operatorname{Ad}_{w^{-1}}(\mathfrak{p}_{\mathcal{I}}^\perp) \)
\[
\pi^{-1}(w) = \frac{f_\mathcal{L} + \mathfrak{b} \cap \operatorname{Ad}_{w^{-1}}(\mathfrak{p}_{\mathcal{I}}^\perp)}{U \cap w^{-1} P_\mathcal{I} w}.
\]

\subsubsection{$T$-fixed points of a particular example}
\label{Section $T$-fixed points of a particular example}
Here is a standard fact used in the proof of the next proposition:

\begin{lemma}\label{standard fact}
    Let \( A \) be a unipotent algebraic group and let \( X \) be an affine \( A \)-variety. Then every \( A \)-orbit in \( X \) is closed.
\end{lemma}

We can calculate the set-theoretical \( T \)-fixed points of \( Z_{\mathcal{I}, \varnothing} = P_\mathcal{I} \backslash\!\!\backslash T^*G \underset{[\mathfrak{g}^*/G]}{\times} T^*G/\!\!/ U \). 

\begin{proposition}
\label{NegaDimension fixed proposition}
    \( (Z_{\mathcal{I}, \varnothing})^T = W/W_{\mathcal{I}} \)
\end{proposition}

\begin{prf}
    Write
    \[
    Z_{\mathcal{I}, \varnothing} := U \backslash \{(g,x) \in G \times \mathfrak{g} \mid x \in \mathfrak{p}_{\mathcal{I}}^\perp, \operatorname{Ad}_g x \in \mathfrak{b}\}/P_\mathcal{I}.
    \]

    Suppose \( (g,x) \in Z_{\mathcal{I}, \varnothing} \) is a fixed point by \( T \). We can assume \( g = w \in W/W_\mathcal{I} \). For \( t \in T \), there exists \( u \in U \) and \( p \in P_\mathcal{I} \) such that \( tw = uwp \) and \( x = \operatorname{Ad}_{p^{-1}} x \). Let \( y = \operatorname{Ad}_w x \); then it leads to \( \operatorname{Ad}_u y = \operatorname{Ad}_t y \). Since \( y \in \operatorname{Ad}_w \mathfrak{p}_{\mathcal{I}}^\perp \cap \mathfrak{b} \), we have \( y \in \mathfrak{u} \).

    Now pick \( h \in \mathfrak{t} \) such that \( [h, e_i] = 2 e_i \). The element \( h \) induces a one-parameter subgroup \( \rho: \mathbb{C}^\times \rightarrow T, z \mapsto \rho(z) \) such that \( \operatorname{Ad}_{\rho(z)} e_{\alpha} = z^{\alpha(h)} e_\alpha \) for \( e_\alpha \in \mathfrak{g}_\alpha \).

    For each \( z \in \mathbb{C}^\times \), there exists \( u(z) \in U \) such that \( \operatorname{Ad}_{\rho(z)} y = \operatorname{Ad}_{u(z)} y \). Since \( y \in \mathfrak{u} \), we have \( \lim_{z \to 0} \operatorname{Ad}_{\rho(z)} y = 0 \). This means the \( U \)-orbit containing \( y \) has a limiting point \( 0 \). By Lemma \ref{standard fact}, every \( U \)-orbit in an affine \( U \)-variety is closed, so \( y = 0 \).
\end{prf}

But it is hard to generalize to calculate the \( T_{f_{\mathcal{L}}} \)-fixed points of \( Z_{\mathcal{I}, \mathcal{L}} \), where \( T_{f_{\mathcal{L}}} := \{t \in T \mid \operatorname{Ad}_t f_\mathcal{L} = f_\mathcal{L}\} = \{t \in T \mid \alpha_j(t) = 1, j \in \mathcal{L}\} \).

We may fix this issue by introducing the \( \mathbb{C}^\times \)-action: \( z \cdot (g,x) = (\rho(z)g, z^2 x) \), \( z \in \mathbb{C}^\times \) for \( (g,x) \in U \backslash \{(g,x) \in G \times \mathfrak{g} \mid x \in \mathfrak{p}_{\mathcal{I}}^\perp, \operatorname{Ad}_g x \in f_\mathcal{L} + \mathfrak{b}\}/P_\mathcal{I} \). Since \( \operatorname{Ad}_{\rho(z)} f_\mathcal{L} = z^{-2} f_\mathcal{L} \), this action is well-defined.

However, some aspects of this problem are still unclear, and we invite further discussion to explore potential solutions.

\subsection{Intersection of two quiver gauge theories}\label{Section intersection quiver}
In this subsection, we shall define the notion of intersection of quiver gauge theories and explain how to view the Higgs branch of the star-shaped quiver gauge theories studied in \cite{gannon2023differential} and \cite{gannon2025functoriality} as intersections of twisted cotangent bundles.

\subsubsection{Quiver gauge theories}
\label{section quiver gauge theories}
Let \( Q \) be a quiver with vertex set \( Q_0 \) and edge set \( Q_1 \). Fix two dimension vectors \( \mathbf{n}: Q_0 \rightarrow \mathbb{Z}_{\geq 0} \) and \( \mathbf{m}: Q_0 \rightarrow \mathbb{Z}_{\geq 0} \). Write \( n_i := \mathbf{n}(i) \) and \( m_i := \mathbf{m}(i) \). In this case, we can consider the usual quiver gauge theory \cite{Braverman:2016pwk}. That is, we let 
\[
V_Q := \bigoplus_{i \in Q_0} \mathbb{C}^{n_i}, \quad W_Q := \bigoplus_{i \in Q_0} \mathbb{C}^{m_i},
\]
which we view as \( Q_0 \)-graded vector spaces. We call \( (V_Q, W_Q) \) a quiver \textit{gauge theory datum} of type \( Q \) and also denote it by \( Q \).

Let 
\[
G(Q) := \mathrm{GL}(V_Q) = \prod_{i \in Q_0} \mathrm{GL}(V_i).
\]

Let 
\[
\mathbf{N}_Q := \bigoplus_{(i \rightarrow j) \in Q_1} \operatorname{Hom}(\mathbb{C}^{n_i}, \mathbb{C}^{n_j}) \oplus \bigoplus_{i \in Q_0} \operatorname{Hom}(\mathbb{C}^{n_i}, \mathbb{C}^{m_i}).
\]

Then \( \mathrm{GL}(V_Q) \) acts on \( T^*\mathbf{N}_Q \) with moment map \( \mu: T^*\mathbf{N}_Q \rightarrow \bigoplus_{i \in Q_0} \mathfrak{gl}_{n_i}(\mathbb{C}) \). The quiver variety \( \mathfrak{M}(V_Q, W_Q) \) (or \( \mathfrak{M}(Q) \) for short) is defined to be the quotient \( \mu^{-1}(0)^{ss}/\mathrm{GL}(V_Q) \), where \( \mu^{-1}(0)^{ss} \subset \mu^{-1}(0) \) consists of those semistable points under some stability condition (for details of the definition, see for example \cite{nakajima2016introduction} and the textbook \cite{kirillov2016quiver}). The quiver variety \( \mathfrak{M}(Q) \) is the Higgs branch of this quiver gauge theory.

Set \( M_C(Q) = M_C(\mathrm{GL}(V_Q), \mathbf{N}_Q) \) to be the Coulomb branch of this quiver gauge theory. 

Let \( \bar{v} = (v_{n+1}, v_n, \ldots, v_1) \) be a vector of decreasing sequences of positive integers. Let \( \mathcal{F}_{\bar{v}} \) denote the partial flag variety parametrizing nested subspaces 
\[
0 \subset \mathcal{V}_1^{v_1} \subset \mathcal{V}_2^{v_2} \subset \ldots \subset \mathcal{V}_{n-1}^{v_{n-1}} \subset \mathcal{V}_n^{v_n} \subset \mathbb{C}^{v_{n+1}},
\]
where \( \mathcal{V}_i^{v_i} \) is a \( v_i \)-dimensional linear subspace of \( \mathbb{C}^{v_{n+1}} \).

Let \( P_{\bar{v}} \) be the parabolic subgroup of \( \operatorname{GL}_{v_{n+1}} \) that fixes the nested subspaces generated by the standard basis of \( \mathbb{C}^{v_{n+1}} \). Then we have the isomorphism \( \mathcal{F}_{\bar{v}} = \operatorname{GL}_{v_{n+1}} / P_{\bar{v}} \).

It is easy to verify that \( T^*\mathcal{F}_{\bar{v}} \) is the same as the quiver variety\footnote{Since \( T^*\mathbf{N}_Q \) doubles the space \( \mathbf{N}_Q \), the direction of the arrows in a quiver does not matter, and we shall ignore this issue in the picture here and below.} associated with 
\[
Q_{\bar{v}} = 
\begin{tikzpicture}[baseline=-20, scale=.4]
 \draw[thick] (1.5,-1) -- (3.5,-1);  
 \draw[thick] (5.6,-1) -- (9,-1); 
  \draw[fill] (1.5,-1) circle (3pt);  
 \draw[fill] (3,-1) circle (3pt);  
    \node at (4.5,-1) {$\cdots$};  
 \draw[fill] (6,-1) circle (3pt); 
 \draw[fill] (7.5,-1) circle (3pt);  
 \draw[fill] (9,-1) circle (3pt);  

 \draw[thick] (8.9,-2.1) rectangle (9.1,-1.9);
 \draw[thick] (9,-1) -- (9,-1.9);
 \node at (1.5,-.5) {$\scriptscriptstyle v_1$} ;
 \node at (3,-.5) {$\scriptscriptstyle v_2$} ;
 \node at (6,-.5) {$\scriptscriptstyle v_{n-2}$} ;
 \node at (7.5,-.5) {$\scriptscriptstyle v_{n-1}$} ;
 \node at (9,-.5) {$\scriptscriptstyle v_{n}$} ;
 \node at (9,-2.6) {$\scriptscriptstyle v_{n+1}$} ;
\end{tikzpicture}.
\]

Here \( Q_{\bar{v}} \) has \( n \) vertices with \( n_i = v_i \) for \( i = 1, \ldots, n \), \( m_i = 0 \) for \( i = 1, \ldots, n-1 \), and \( m_n = v_{n+1} \).

\subsubsection{Intersection of two quiver gauge theories}
Suppose we have two quiver gauge theories, \( Q \) and \( Q' \), that have framing \( m_i \) and \( m_{i'} \) at the vertices \( i \) and \( i' \). Say 
\[
Q = \begin{tikzpicture}[baseline=-20, scale=.4]
 \draw[thick] (1.5,-1.5) -- (4,-1.5);  
 \draw[thick] (3,-1.5) -- (4,-1); 
 \draw[thick] (3,-1.5) -- (4,-2); 
 \draw[thick] (1.3,-1.6) rectangle (1.5,-1.4);
 \draw[fill] (3,-1.5) circle (3pt);  
    \node at (5,-1.5) {$\cdots$};  
 \node at (5,-1) {$\cdots$};  
 \node at (5,-2) {$\cdots$};  
 \node at (1.5,-1) {$\scriptscriptstyle m_i$} ;
 \node at (3,-1) {$\scriptscriptstyle n_i$} ;
\end{tikzpicture}
\]
and 
\[
Q' = \begin{tikzpicture}[baseline=-20, scale=.4]
\draw[thick] (1.3,-1.6) rectangle (1.5,-1.4);
 \draw[thick] (1.5,-1.5) -- (3,-1.5);  
 \draw[thick] (3,-1.5) -- (4,-1); 
 \draw[thick] (3,-1.5) -- (4,-2); 
 \draw[fill] (3,-1.5) circle (3pt);  
 
 \node at (5,-1) {$\cdots$};  
 \node at (5,-2) {$\cdots$};  
 \node at (1.5,-1) {$\scriptscriptstyle m_{i'}'$} ;
 \node at (3,-1) {$\scriptscriptstyle n_{i'}'$} ;
\end{tikzpicture}
\]
with \( m_i \leq m_{i'}' \).

Then the \textit{intersection} of \( Q \) and \( Q' \) at \( i \) and \( i' \) is defined as the new quiver gauge theory:
\[
Q \times_{(i,i')} Q' := \begin{tikzpicture}[baseline=-20, scale=.4]
   \node at (-1.4,-1.5) {$\cdots$};  
 \node at (-1.4,-1) {$\cdots$};  
 \node at (-1.4,-2) {$\cdots$}; 
 \draw[thick] (-0.8,-1.5) -- (1.5,-1.5);  
  \draw[thick] (-0.8,-2) -- (0,-1.5);  
   \draw[thick] (-0.8,-1) -- (0,-1.5);  
\draw[fill] (1.5,-1.5) circle (3pt);
\draw[fill] (0,-1.5) circle (3pt);
 \draw[thick] (1.5,-1.5) -- (3,-1.5);  
 \draw[thick] (3,-1.5) -- (4,-1); 
 \draw[thick] (3,-1.5) -- (4,-2); 
 \draw[fill] (3,-1.5) circle (3pt);  

\draw[thick] (1.5,-1.5) -- (1.5,-0.5); 
 \draw[thick] (1.4,-0.5) rectangle (1.6,-0.3);
  
 \node at (5,-1) {$\cdots$};  
 \node at (5,-2) {$\cdots$};  
 \node at (0,-1) {$\scriptscriptstyle n_i$} ;
 \node at (1.5,-2) {$\scriptscriptstyle m_i$} ;
 \node at (3,-1) {$\scriptscriptstyle n_{i'}'$} ;

  \node at (1.5,0) {$\scriptscriptstyle m_{i'}' - m_i$} ;
\end{tikzpicture}.
\]

Diagrammatically, we connect the two quivers at the framing \( m_i \), and the remaining framing integer is the absolute value of the difference.

The quiver variety of the intersection \( Q \times_{(i,i')} Q' \) is the same as the intersection of the corresponding quiver varieties \( \mathfrak{M}(Q) \underset{[\mathfrak{gl}_{m_i}^*/\operatorname{GL}_{m_i}]}{\times} \mathfrak{M}(Q') \) if we do not consider any of the stability conditions during the construction of quiver varieties. This is why we call $Q \times_{(i,i')} Q'$ the intersection of $Q$ and $Q'$.

We can also intersect two quiver gauge theories at multiple framing vertices at the same time, and then the intersection of the corresponding quiver varieties is obtained by a diagonal quotient by the product \( \operatorname{GL}_{m_{i_1}} \times \operatorname{GL}_{m_{i_2}} \times \cdots \times \operatorname{GL}_{m_{i_k}} \).

\subsubsection{The star-shaped quiver gauge theories}A graph \( \Gamma \) is \textit{star-shaped} if it is a tree, and there is some vertex \( v \) such that each connected component of \( \Gamma \backslash \{v\} \) is an \( A \)-series Dynkin diagram whose vertex adjacent to \( v \) is one of the two outermost vertices.

We call \( v \) the \textit{center} of \( \Gamma \). Each connected component of \( \Gamma \backslash \{v\} \) is referred to as a \textit{leg} of \( \Gamma \). The center is unique if and only if the number of legs is no less than 3.

A quiver gauge theory is star-shaped if its underlying graph is star-shaped. We are mostly interested in the star-shaped quiver gauge theories discussed in \cite{gannon2025functoriality} and \cite{gannon2023differential}.

Let \( \bar{n} = (n_1, n_2, \ldots, n_k) \) be an ordered tuple of positive integers summing to \( n \). Let 
\[
Q_{\bar{n}}^{\star} := \begin{tikzpicture}[baseline=-20, scale=.4]
   \node at (-0.5,-1.5) {$\cdots$};  
 \draw[thick] (-3,-1.5) -- (-1.2,-1.5);  
 \draw[thick] (0.3,-1.5) -- (3,-1.5);  
 \draw[thick] (3,-1.5) -- (5,0); 
 \draw[thick] (3,-1.5) -- (5,-0.5); 
  \draw[thick] (3,-1.5) -- (5,-3); 
  \draw[fill] (-3,-1.5) circle (3pt);
    \draw[fill] (-1.5,-1.5) circle (3pt);
 \draw[fill] (1.5,-1.5) circle (3pt);

 \draw[fill] (3,-1.5) circle (3pt);  
  \draw[fill] (5,0) circle (3pt);  
   \draw[fill] (5,-0.5) circle (3pt); 
   \draw[fill] (5,-3) circle (3pt);  
 
 \node at (5,-1.5) {$\vdots$};  
 \node at (-3,-1) {$\scriptscriptstyle 1$} ;
 \node at (-1.5,-1) {$\scriptscriptstyle 2$} ;
 \node at (1.5,-1) {$\scriptscriptstyle n-2$} ;
 \node at (3,-1) {$\scriptscriptstyle n-1$} ;
 \node at (5.9,0) {$\scriptscriptstyle n_1$} ;
 \node at (5.9,-0.5) {$\scriptscriptstyle n_2$} ;
 \node at (5.9,-3) {$\scriptscriptstyle n_k$} ;
\end{tikzpicture}.
\]Then in \cite{gannon2023differential}, it is proved that 
\[
M_C(Q_{\bar{n}}^{\star}) = T^*\operatorname{GL}_n(\mathbb{C})/\!\!/_{\psi_{\bar{n}}} U,
\]
where \( \psi_{\bar{n}}: \mathfrak{u} \rightarrow \mathbb{C} \) is the character such that \( \psi_{\bar{n}}(e_i) = 0 \) if \( i = n_1 + \cdots + n_j \) for some \( j \), or \( \psi_{\bar{n}}(e_i) = 1 \) otherwise.

Let $
L_{\bar{n}} := \operatorname{GL}_{n_1} \times \operatorname{GL}_{n_2} \times \cdots \times \operatorname{GL}_{n_k}$ be the Levi subgroup of \( \operatorname{GL}_n \). Then we also view \( \mathfrak{M}(Q_{\bar{n}}^{\star}) \) as the intersection
\[
\mathfrak{M}(Q_{\bar{n}}^{\star}) = \mathfrak{M}(Q_{(n,n-1,\ldots,1)}) \underset{[\mathfrak{l}_{\bar{n}}^*/{L}_{\bar{n}}]}{\times} \mathfrak{M} \left( \begin{tikzpicture}[baseline=-20, scale=.4]
   
  \draw[fill] (1,0) circle (3pt);  
   \draw[fill] (1,-0.5) circle (3pt); 
   \draw[fill] (1,-3) circle (3pt);  
   \draw[thick] (1,0) -- (2.5,0); 
   \draw[thick] (1,-0.5) -- (2.5,-0.5); 
   \draw[thick] (1,-3) -- (2.5,-3); 
   \draw[thick] (2.5,-0.1) rectangle (2.7,0.1);
   \draw[thick] (2.5,-0.6) rectangle (2.7,-0.4);
   \draw[thick] (2.5,-3.1) rectangle (2.7,-2.9);
   \node at (1,-1.5) {$\vdots$};  
   \node at (2.5,-1.5) {$\vdots$};  
   \node at (0.3,0) {$\scriptscriptstyle 0$} ;
   \node at (0.3,-0.5) {$\scriptscriptstyle 0$} ;
   \node at (0.3,-3) {$\scriptscriptstyle 0$} ;
   \node at (3.5,0) {$\scriptscriptstyle n_1$} ;
   \node at (3.5,-0.5) {$\scriptscriptstyle n_2$} ;
   \node at (3.5,-3) {$\scriptscriptstyle n_k$} ;
\end{tikzpicture} \right).
\]

Now we have 
\[
\mathfrak{M}(Q_{(n,n-1,\ldots,1)}) = B \backslash\!\!\backslash T^*\operatorname{GL}_n
\]
and 
\[
\mathfrak{M} \left( \begin{tikzpicture}[baseline=-20, scale=.4]
   
  \draw[fill] (1,0) circle (3pt);  
   \draw[fill] (1,-0.5) circle (3pt); 
   \draw[fill] (1,-3) circle (3pt);  
   \draw[thick] (1,0) -- (2.5,0); 
   \draw[thick] (1,-0.5) -- (2.5,-0.5); 
   \draw[thick] (1,-3) -- (2.5,-3); 
   \draw[thick] (2.5,-0.1) rectangle (2.7,0.1);
   \draw[thick] (2.5,-0.6) rectangle (2.7,-0.4);
   \draw[thick] (2.5,-3.1) rectangle (2.7,-2.9);
   \node at (1,-1.5) {$\vdots$};  
   \node at (2.5,-1.5) {$\vdots$};  
   \node at (0.3,0) {$\scriptscriptstyle 0$} ;
   \node at (0.3,-0.5) {$\scriptscriptstyle 0$} ;
   \node at (0.3,-3) {$\scriptscriptstyle 0$} ;
   \node at (3.5,0) {$\scriptscriptstyle n_1$} ;
   \node at (3.5,-0.5) {$\scriptscriptstyle n_2$} ;
   \node at (3.5,-3) {$\scriptscriptstyle n_k$} ;
\end{tikzpicture} \right) = \mathfrak{M} \left( \coprod_{i=1}^k Q_{(n_i)} \right) = T^*L_{\bar{n}}/\!\!/L_{\bar{n}}.
\]
Therefore,
\[
\mathfrak{M}(Q_{\bar{n}}^{\star}) = B \backslash\!\!\backslash T^*\operatorname{GL}_n \underset{[\mathfrak{l}_{\bar{n}}^*/{L}_{\bar{n}}]}{\times} T^*L_{\bar{n}}/\!\!/ L_{\bar{n}} = B \backslash\!\!\backslash T^*\operatorname{GL}_n \underset{[\mathfrak{gl}_{{n}}^*/\operatorname{GL}_n]}{\times}T^*\operatorname{GL}_n/\!\!/ L_{\bar{n}}.
\]

Now we consider another star-shaped quiver gauge theory studied in \cite{gannon2025functoriality}:
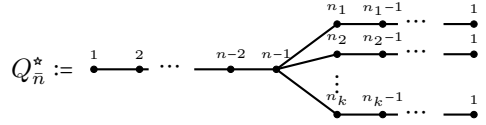
\[
Q^\smallstar_{\bar{n}} := \begin{tikzpicture}[baseline=-20, scale=.4]
   \node at (-0.5,-1.5) {$\cdots$};  
 \draw[thick] (-3,-1.5) -- (-1.2,-1.5);  
 \draw[thick] (0.3,-1.5) -- (3,-1.5);  
 \draw[thick] (3,-1.5) -- (5,0); 
 \draw[thick] (3,-1.5) -- (5,-1); 
 \draw[thick] (3,-1.5) -- (5,-3); 

 \draw[thick] (5,0) -- (7,0); 
 \draw[thick] (8.5,0) -- (9.5,0); 
 \node at (7.6,0) {$\cdots$};  
 \draw[fill] (6.5,0) circle (3pt); 
 \draw[fill] (9.5,0) circle (3pt); 
 \node at (5,0.5) {$\scriptscriptstyle n_1$} ;
 \node at (6.5,0.5) {$\scriptscriptstyle n_1-1$} ;
 \node at (9.5,0.5) {$\scriptscriptstyle 1$} ;

 \draw[thick] (5,-1) -- (7,-1); 
 \draw[thick] (8.5,-1) -- (9.5,-1); 
 \node at (7.6,-1) {$\cdots$};  
 \draw[fill] (6.5,-1) circle (3pt); 
 \draw[fill] (9.5,-1) circle (3pt);
 \node at (5,-0.5) {$\scriptscriptstyle n_2$} ;
 \node at (6.5,-0.5) {$\scriptscriptstyle n_2-1$} ;
 \node at (9.5,-0.5) {$\scriptscriptstyle 1$} ;

 \draw[thick] (5,-3) -- (7,-3); 
 \draw[thick] (8.5,-3) -- (9.5,-3); 
 \node at (7.6,-3) {$\cdots$};  
 \draw[fill] (6.5,-3) circle (3pt); 
 \draw[fill] (9.5,-3) circle (3pt);
 \node at (5,-2.5) {$\scriptscriptstyle n_k$} ;
 \node at (6.5,-2.5) {$\scriptscriptstyle n_k-1$} ;
 \node at (9.5,-2.5) {$\scriptscriptstyle 1$} ;

 \draw[fill] (-3,-1.5) circle (3pt);
 \draw[fill] (-1.5,-1.5) circle (3pt);
 \draw[fill] (1.5,-1.5) circle (3pt);
 \draw[fill] (3,-1.5) circle (3pt);  
 \draw[fill] (5,0) circle (3pt);  
 \draw[fill] (5,-1) circle (3pt); 
 \draw[fill] (5,-3) circle (3pt);  

 \node at (5,-2) {$\vdots$};  
 \node at (-3,-1) {$\scriptscriptstyle 1$} ;
 \node at (-1.5,-1) {$\scriptscriptstyle 2$} ;
 \node at (1.5,-1) {$\scriptscriptstyle n-2$} ;
 \node at (3,-1) {$\scriptscriptstyle n-1$} ;
\end{tikzpicture}.
\]Its Coulomb branch is the affinization of \( T^*\operatorname{GL}_n/\!\!/U_{\bar{n}} \). The unipotent subgroup \( U_{\bar{n}} \) is the unipotent radical  of the group of block upper triangular matrices in \( \operatorname{GL}_n \) with block sizes indexed by the ordered partition \( n = n_1 + \cdots + n_k \) of \( n \).

Similarly, we also have 
\[
\begin{aligned}
    \mathfrak{M}(Q^\smallstar_{\bar{n}}) &= \mathfrak{M}(Q_{(n,n-1,\ldots,1)}) \underset{[\mathfrak{l}_{\bar{n}}^*/{L}_{\bar{n}}]}{\times} \mathfrak{M} \left( \coprod_{i=1}^k Q_{(n_i,n_i-1,\ldots,1)} \right)\\
    &= B \backslash\!\!\backslash T^*\operatorname{GL}_n \underset{[\mathfrak{l}_{\bar{n}}^*/{L}_{\bar{n}}]}{\times} T^*L_{\bar{n}}/\!\!/B_{L_{\bar{n}}}\\
    &= B \backslash\!\!\backslash T^*\operatorname{GL}_n\underset{[\mathfrak{gl}_{{n}}^*/\operatorname{GL}_n]}{\times} T^*\operatorname{GL}_n/\!\!/B_{L_{\bar{n}}}.
\end{aligned}
\]

Thus, \( \mathfrak{M}(Q^\smallstar_{\bar{n}}) \) could also be understood as an intersection of cotangent bundles.

\subsection{Bow varieties as intersections of twisted cotangent bundles}\label{section bow variety}
In this subsection, we shall show that bow varieties with good properties are intersections of twisted cotangent bundles.

Objects like \( \db 2 \nsb 3 \db 1 \db 5 \nsb 7 \db \) are called (type A) brane diagrams and are denoted as \( \mathcal{D} \). A red line $\nsb$ represents a \textit{NS5 brane}, while a blue line $\db$ represents a \textit{D5 brane}. Associated with a brane diagram \( \mathcal{D} \) is a smooth holomorphic symplectic variety \( \mathcal{C}(\mathcal{D}) \) called the bow variety. One definition is through quivers (\cite{nakajima2017cherkis}, \cite{rimanyi2020bow}).

The Hanany-Witten transition is the local modification in a brane diagram
\[
d_1 \db d_2 \nsb d_3 \stackrel{HW}{\longleftrightarrow} d_1 \nsb \tilde{d_2} \db d_3, \quad \tilde{d_2}=d_1+d_3+1-d_2.
\]
which induces an isomorphism between the corresponding bow varieties. We can use the Hanany-Witten transition to move all the D5 branes to the right and NS5 branes to the left. Then typically a brane diagram will look like 
\[
\mathcal{D}_{\hat{m}} = \nsb m_{-u} \nsb \cdots \nsb m_{-1} \nsb m_0 \db m_1 \db \cdots \db m_v \db,
\]
where \( \hat{m} = (m_{-u}, \ldots, m_{-1}, m_0, m_1, \ldots, m_v) \). Let \( \hat{m}^- = (m_0, m_{-1}, \ldots, m_{-u}) \) and \( \hat{m}^+ = (m_0, m_1, \ldots, m_v) \).

If \( 0 < m_{-u} < \cdots < m_{-1} < m_0 > m_1 > \cdots > m_v > 0 \), one particular result proved in \cite{ji2023bow} is that 
\[
\begin{aligned}
    \mathcal{C}(\mathcal{D}_{\hat{m}}) &= \mathcal{C}(\nsb m_{-u} \nsb \cdots \nsb m_{-1} \nsb m_0) \underset{[\mathfrak{gl}_{{m_0}}^*/\operatorname{GL}_{m_0}]}{\times} \mathcal{C}(m_0 \db m_1 \db \cdots \db m_v \db) \\
&= P_{\hat{m}^-} \!\!\backslash\!\!\backslash T^*\operatorname{GL}_{m_0} \underset{[\mathfrak{gl}_{{m_0}}^*/\operatorname{GL}_{m_0}]}{\times}  T^*\operatorname{GL}_{m_0}/\!\!/_{f_{\hat{m}^+}} U_{f_{\hat{m}^+}},
\end{aligned}
\]
where the character \( f_{\hat{m}^+} \) and the unipotent group \( U_{f_{\hat{m}^+}} \) are defined in \cite[Section~3]{ji2023bow}.

By replacing NS5-branes with D5-branes, and D5-branes with NS5-branes in a brane diagram \( \mathcal{D} \), we obtain another brane diagram \( \check{\mathcal{D}} \), called the 3d mirror dual brane diagram to \( \mathcal{D} \). Then  
\[
\begin{aligned}
   \mathcal{C}(\check{\mathcal{D}}_{\hat{m}}) &= \mathcal{C}(\db m_{-u} \db \cdots \db m_{-1} \db m_0) \underset{[\mathfrak{gl}_{{m_0}}^*/\operatorname{GL}_{m_0}]}{\times}  \mathcal{C}(m_0 \nsb m_1 \nsb \cdots \nsb m_v)\\
   &= U_{f_{\hat{m}^-}} {}_{f_{\hat{m}^-}}\!\! \backslash\!\!\backslash T^*\operatorname{GL}_{m_0} \underset{[\mathfrak{gl}_{{m_0}}^*/\operatorname{GL}_{m_0}]}{\times}  T^*\operatorname{GL}_{m_0}/\!\!/ P_{\hat{m}^+}. 
\end{aligned}
\]

Now we turn to the extremal cases, where a bow diagram consists of either D5 branes or NS5 branes.

Let \( D^{\diamond}_{\hat{m}} = \db m_{-u} \db \cdots \db m_{-1} \db m_0 \db m_1 \db \cdots \db m_v \db \). Then if we don't delete the points that are unstable at $m_0$ in the sense of \cite[p.~10]{rimanyi2020bow} when taking quotients, we have 
\[
\begin{aligned}
    \mathcal{C}(D^{\diamond}_{\hat{m}}) &= \mathcal{C}(\db m_{-u} \db \cdots \db m_{-1} \db m_0) \underset{[\mathfrak{gl}_{{m_0}}^*/\operatorname{GL}_{m_0}]}{\times}  \mathcal{C}(m_0 \db m_1 \db \cdots \db m_v) \\
    &= U_{f_{\hat{m}^-}} {}_{f_{\hat{m}^-}} \!\!\backslash\!\!\backslash T^*\operatorname{GL}_{m_0} \underset{[\mathfrak{gl}_{{m_0}}^*/\operatorname{GL}_{m_0}]}{\times}  T^*\operatorname{GL}_{m_0}/\!\!/_{f_{\hat{m}^+}} U_{f_{\hat{m}^+}}.
\end{aligned}
\]

The 3d mirror variety to $\mathcal{C}(D^{\diamond}_{\hat{m}})$ is the bow variety associated with the 3d mirror dual brane diagram $\check{D}^{\diamond}_{\hat{m}}$:
\[
\begin{aligned}
    \mathcal{C}(\check{D}^{\diamond}_{\hat{m}}) &= \mathcal{C}(\nsb m_{-u} \nsb \cdots \nsb m_{-1} \nsb m_0) \underset{[\mathfrak{gl}_{{m_0}}^*/\operatorname{GL}_{m_0}]}{\times}  \mathcal{C}(m_0 \nsb m_1 \nsb \cdots \nsb m_v) \\
    &= P_{\hat{m}^-}\!\! \backslash\!\!\backslash T^*\operatorname{GL}_{m_0} \underset{[\mathfrak{gl}_{{m_0}}^*/\operatorname{GL}_{m_0}]}{\times}  T^*\operatorname{GL}_{m_0}/\!\!/ P_{\hat{m}^+}.
\end{aligned}
\]
$\mathcal{C}(\check{D}^{\diamond}_{\hat{m}})$ can also be viewed as the quiver variety of the self-intersection of the quiver $Q_{\bar{v}}$ defined in Section \ref{section quiver gauge theories}, which is also a star-shaped quiver.

\begin{example}
    If \( \hat{m} = (m_0) \), then in this case it is a reformulation of Proposition \ref{UniversalCent as intersection} in the language of bow varieties that 
\[
\mathcal{C}(\db m_0 \db) = \{(g,x) \in \operatorname{GL}_{m_0}(\mathbb{C}) \times \mathfrak{gl}_{m_0}(\mathbb{C}) \mid \operatorname{Ad}_g x = x\}
\]
is the universal centralizer of $\operatorname{GL}_{m_0}$. Its 3d mirror variety is 
\[
\operatorname{GL}_{m_0}\!\! \backslash\!\!\backslash T^*\operatorname{GL}_{m_0} \underset{[\mathfrak{gl}_{{m_0}}^*/\operatorname{GL}_{m_0}]}{\times}  T^*\operatorname{GL}_{m_0}/\!\!/ \operatorname{GL}_{m_0}.
\]
\end{example}

\section{Closed subsets of root systems}
In this section, we study the combinatorial aspect of subgroups with maximal rank in \( G \).

Recall from \cite{bourbaki1994lie}, a subset \( \Gamma \subset \Phi \) is \textit{closed} if \( \alpha, \beta \in \Gamma \) and \( \alpha + \beta \in \Phi \) implies \( \alpha + \beta \in \Gamma \). A subset \( \Gamma \) is closed if and only if \( -\Gamma \) is closed. Closed subsets of $\Phi$ correspond to connected subgroups \( H \subset G \) such that \( T \subset H \). We will write the corresponding connected subgroup \( H \) as \( H_\Gamma \). The Lie algebra of \( H_\Gamma \) is \( \mathfrak{h}_\Gamma = \mathfrak{t} \oplus_{\alpha \in \Gamma} \mathfrak{g}_\alpha \). Then its perpendicular complement with respect to the killing form is
\[
\mathfrak{h}_\Gamma^\perp = \mathfrak{t}^\perp \cap \bigcap_{\alpha \in \Gamma} \mathfrak{g}_\alpha^\perp = \mathfrak{t}^\perp \cap \bigcap_{\alpha \in \Gamma} \left( \mathfrak{t} \oplus \bigoplus_{\beta \in \Phi \backslash \{-\alpha\}} \mathfrak{g}_\beta \right) = \bigoplus_{\Phi \backslash \{-\Gamma\}} \mathfrak{g}_\alpha = V_{\Phi \backslash \{-\Gamma\}}.
\]

The intersection of closed subsets is closed, and the corresponding subgroup is the identity component of the intersection of the corresponding subgroups.

A closed subset \( \Gamma \) is \textit{invertible} if \( \Phi \backslash \Gamma \) is closed.

A subset \( J \) of \( \Gamma \) is called an \textit{ideal} if 
\[
\alpha \in \Gamma, \ \beta \in J, \ \alpha + \beta \in \Gamma \Rightarrow \alpha + \beta \in J.
\]

Let \( \Gamma_s := \{ \alpha \in \Gamma \mid -\alpha \in \Gamma \} \) and \( \Gamma_u = \Gamma \backslash \Gamma_s \).

It is easy to see that \( \Gamma_s \) is closed and \( \Gamma_u \) is an ideal of \( \Gamma \). We will call \( \Gamma_s \) the symmetric part of \( \Gamma \) and \( \Gamma_u \) the unipotent part of \( \Gamma \). If \( \Gamma = \Gamma_s \), then \( \Gamma \) is called \textit{symmetric}. If \( \Gamma = \Gamma_u \), then \( \Gamma \) is called \textit{unipotent}.

\begin{definition}
    Closed subsets \( \Gamma \) and \( \Gamma' \) are conjugate if there exists an element \( w \in W \) such that \( w(\Gamma) = \Gamma' \).
\end{definition}

Recall \cite{bourbaki1994lie}, Chapter VI, Proposition 22: the unipotent part \( \Gamma_u \) is conjugate to a subset of \( \Phi^+ \). Let \( \Gamma \) be a unipotent closed subset of \( \Phi \). We will denote \( U_{\Gamma} \) the unipotent group corresponding to the nilpotent Lie algebra \( \mathfrak{u}_\Gamma = \oplus_{\alpha \in \Gamma} \mathfrak{g}_\alpha \).

If \( \Gamma \cup (-\Gamma) = \Phi \), we will call \( \Gamma \) a \textit{parabolic} subset.

Let \( \mathcal{I} \subset \Pi \). Then a parabolic subset is conjugate to a closed subset of the form \( \Phi_\mathcal{I} \cup \Phi^+ \) (see \cite{bourbaki1994lie}, Chapter VI, Proposition 21). However, the \( \Phi_\mathcal{I} \)'s are not all representatives of \( W \)-orbits of symmetric subsets of \( \Phi \) (for example, \( A_1 \times A_1 \) in \( B_2 \) is not conjugate to the standard Levi's).

\subsection{A family of closed subsets}

Let \( (\mathcal{I}, \mathcal{J}, \mathcal{K}) \) be a triple of subsets of \( \Pi \) such that \( \mathcal{I}, \mathcal{K} \subset \mathcal{J} \). Define 
\[
\Gamma(\mathcal{I}, \mathcal{J}, \mathcal{K}) := \Phi_{\mathcal{I}} \cup (\Phi^+_\mathcal{J} \backslash \Phi_\mathcal{K}).
\]

When \( \mathcal{J} = \Pi \), we will abbreviate \( \Gamma(\mathcal{I}, \Pi, \mathcal{K}) \) as \( \Gamma(\mathcal{I}, \mathcal{K}) \).

When \( \mathcal{I} = \mathcal{K} \) and \( \mathcal{J} = \Pi \), \( P(\mathcal{I}) := \Gamma(\mathcal{I}, \mathcal{I}) = \Phi_{\mathcal{I}} \cup \Phi^+ \) is parabolic. When \( \mathcal{K} = \Pi \), \( \Gamma(\mathcal{I}, \Pi) = \Phi_{\mathcal{I}} \). When \( \mathcal{I} = \mathcal{K} \), \( \Gamma(\mathcal{I}, \mathcal{J}, \mathcal{I}) \) is a parabolic subset of the subsystem \( \Phi_\mathcal{J} \).

Therefore, \( H_{\Gamma(\mathcal{I}, \mathcal{J}, \mathcal{K})} \) could be viewed as a discontinuous deformation between parabolic subgroups of Levi subgroups. All the subgroups of \( G \) that contain \( T \) and are considered in Section \ref{Intersection of Hamiltonian $G$-varieties} are of this form.

We classify when \( \Gamma(\mathcal{I}, \mathcal{J}, \mathcal{K}) \) is closed:

\begin{proposition}\label{general closed condition}
    \( \Gamma(\mathcal{I}, \mathcal{J}, \mathcal{K}) \) is closed if and only if there exist \( \mathcal{X}, \mathcal{Y} \subset \mathcal{J} \) such that the following are satisfied:
    \begin{enumerate}
        \item \( \mathcal{I} \cup \mathcal{K} = \mathcal{X} \coprod \mathcal{Y} \)
        \item \( \mathcal{X} \) is orthogonal to \( \mathcal{Y} \),
        \item \( \mathcal{X} \supset \mathcal{I} \backslash \mathcal{K} \), \( \mathcal{Y} \supset \mathcal{K} \backslash \mathcal{I} \).
    \end{enumerate}
\end{proposition}

\begin{prf}We may assume $\mathcal{J}=\Pi$ by replacing $G$ with $H_{\Phi_\mathcal{J}}$.
    
    Suppose $\mathcal{X},\mathcal{Y}$ are subsets of $\Pi$ satisfying the three conditions above. Let $\alpha,\beta\in \Gamma(\mathcal{I},\mathcal{K})$ such that $\alpha+\beta\in\Phi$. Since $\Phi_\mathcal{I}^+\cup\Phi^+\backslash\Phi_\mathcal{K}$ is closed, we can assume $\alpha\in\Phi^-_\mathcal{I},\beta\in \Phi^+\backslash\Phi_\mathcal{K}$. Write $\beta=\beta_1+\beta_2+\beta_3+\beta_4,$ where $ \beta_1\in\langle \mathcal{I}\backslash \mathcal{K}\rangle,\beta_2\in\langle \mathcal{I}\cap \mathcal{K}\rangle,\beta_3\in\langle \mathcal{K}\backslash \mathcal{I}\rangle,\beta_4\in\langle \Pi\backslash(\mathcal{I}\cup \mathcal{K})\rangle$ are unique. If $\beta_4\neq 0$, then $\alpha+\beta \in \Phi^+\backslash \Phi_\mathcal{K}\subset \Gamma(\mathcal{I},\mathcal{K})$. If $\beta_4=0$, then $\beta_1\neq 0$. Since now $\beta\in \Phi_{\mathcal{I}\cup \mathcal{K}}$ and $\Phi_{\mathcal{I}\cup \mathcal{K}}=\Phi_\mathcal{X}\coprod\Phi_\mathcal{Y}$, we must have $\beta_3=0$, and thus $\alpha+\beta\in \Phi\cap\langle \mathcal{I}\rangle=\Phi_\mathcal{I}\subset \Gamma(\mathcal{I},\mathcal{K})$.

    Conversely, suppose $\Gamma(\mathcal{I},\mathcal{K})$ is closed. If $\mathcal{I}\backslash\mathcal{K}$ is not orthogonal to $\mathcal{K}\backslash\mathcal{I}$, then we could pick $\alpha\in \mathcal{I}\backslash\mathcal{K}$ and $\beta\in\mathcal{K}\backslash\mathcal{I}$ such that $(\alpha,\beta)<0$. Then $\alpha+\beta$ is a root and lies in $\Phi^+\backslash\Phi_\mathcal{K}$. But $\beta=(-\alpha)+(\alpha+\beta)$ is a root and is not in $\Gamma(\mathcal{I},\mathcal{K})$, which is a contradiction. So we have $\mathcal{I}\backslash\mathcal{K}$ is orthogonal to $\mathcal{K}\backslash\mathcal{I}$.
Let $\mathcal{I}_0=\mathcal{I}$, $\mathcal{K}_0=\mathcal{K}$.
After defining $\mathcal{I}_k,\mathcal{K}_k$, let
$$\mathcal{I}_{k+1}:=\{\alpha\in \mathcal{I}_k\cap\mathcal{K}_k\mid \alpha \text{ is orthogonal to }\mathcal{K}_k\backslash\mathcal{I}_k\}$$
$$\mathcal{K}_{k+1}:=\{\alpha\in \mathcal{I}_k\cap\mathcal{K}_k\mid \alpha \text{ is orthogonal to }\mathcal{I}_k\backslash\mathcal{K}_k\}.$$
Assume $\mathcal{I}_{k}\backslash\mathcal{K}_{k}$ is orthogonal to $\mathcal{K}_{k}\backslash\mathcal{I}_{k}$. We claim that $\mathcal{I}_{k+1}\backslash\mathcal{K}_{k+1}$ is orthogonal to $\mathcal{K}_{k+1}\backslash\mathcal{I}_{k+1}$.  In fact, if $\alpha_{k+1}\in\mathcal{I}_{k+1}\backslash\mathcal{K}_{k+1},\beta_{k+1}\in\mathcal{K}_{k+1}\backslash\mathcal{I}_{k+1}$ and $(\alpha_{k+1},\beta_{k+1})\neq 0$, we have $\alpha_{k+1}+\beta_{k+1}\in\Phi$. And we could choose $\alpha_k\in \mathcal{I}_{k}\backslash\mathcal{K}_{k}$ such that $(\alpha_k,\alpha_{k+1})\neq 0$, and $\beta_k\in \mathcal{K}_{k}\backslash\mathcal{I}_{k}$ such that $(\beta_k,\beta_{k+1})\neq 0$. Then $\alpha_{k}+\alpha_{k+1}+\beta_{k+1}+\beta_k\in \Phi$. Inductively, we can find $\alpha_{k-1},\alpha_{k-1},\cdots,\alpha_0$ and $\beta_{k-1},\beta_{k-1},\cdots,\beta_0$, such that $\alpha_i\in\mathcal{I}_{i}\backslash\mathcal{K}_{i}$,  $\beta_i\in\mathcal{K}_{i}\backslash\mathcal{I}_{i}$ and $\alpha_0+\alpha_1+\cdots+\alpha_{k+1}+\beta_{k+1}+\cdots+\beta_1, \alpha_0+\alpha_1+\cdots+\alpha_{k+1}+\beta_{k+1}+\cdots+\beta_0\in\Phi$. Now $-(\alpha_0+\alpha_1+\cdots+\alpha_{k+1}+\beta_{k+1}+\cdots+\beta_1)\in \Phi_\mathcal{I}$,  $\alpha_0+\alpha_1+\cdots+\alpha_{k+1}+\beta_{k+1}+\cdots+\beta_0\in \Phi^+\backslash\Phi_\mathcal{K}$, this would imply $\beta_0\in\Gamma(\mathcal{I},\mathcal{K})$, which is a contradiction.

Let 
$$\mathcal{X}=\coprod_{i=0}^{\infty}\mathcal{I}_i\backslash\mathcal{K}_i,\ \ \mathcal{Y}=\coprod_{i=0}^{\infty}\mathcal{K}_i\backslash\mathcal{I}_i.$$

Since $\mathcal{I}_{k+1}\backslash\mathcal{K}_{k+1}$ is also orthogonal to $\mathcal{K}_{i+1}\backslash\mathcal{I}_{i+1}$ for $i< k$ by definition, we have $\mathcal{X},\mathcal{Y}$ defined above satisfying 1. $\mathcal{I}\cup\mathcal{K}=\mathcal{X}\coprod \mathcal{Y}$, 2. $\mathcal{X}$ is orthogonal to $\mathcal{Y}$, 3. $\mathcal{X}\supset \mathcal{I}\backslash \mathcal{K}$, $\mathcal{Y}\supset \mathcal{K}\backslash \mathcal{I}$.
\end{prf}
\begin{problem}
    What is the condition about $(\mathcal{I}, \mathcal{J}, \mathcal{K})$ for $\Gamma(\mathcal{I}, \mathcal{J}, \mathcal{K})$ to be closed when $\mathcal{I} \not\subset \mathcal{J}$?
\end{problem}

The definition of $\Gamma(\mathcal{I}, \mathcal{J}, \mathcal{K})$ is motivated by \cite{dokovic1994closed}, where they defined $\Gamma(\mathcal{I}, \mathcal{K})$ for $\mathcal{I} \subset \mathcal{K}$ and when $\mathcal{I}$ is orthogonal to $\mathcal{K} \backslash \mathcal{I}$. They show such subsets are representatives of invertible subsets under conjugation. This result is proved in \cite[Theorem~4]{dokovic1994closed}.

\begin{theorem}
    Every invertible subset of $\Phi$ is conjugate to $\Gamma(\mathcal{I}, \mathcal{K})$ for some subsets $\mathcal{I} \subset \mathcal{K} \subset \Pi$ and $\mathcal{I}$ orthogonal to $\mathcal{K} \backslash \mathcal{I}$.
\end{theorem}

We prove that there are no other invertible subsets in the family of closed subsets $\Gamma(\mathcal{I}, \mathcal{K})$ for $\mathcal{I}, \mathcal{K} \subset \Pi$.

\begin{proposition}
    $\Gamma(\mathcal{I}, \mathcal{K})$ is invertible if and only if $\mathcal{I}$ is orthogonal to $\mathcal{K} \backslash \mathcal{I}$.
\end{proposition}

\begin{prf}
    We have 
    $$\Phi \backslash \Gamma(\mathcal{I}, \mathcal{K}) = \Phi \backslash \Phi_\mathcal{I} \cap (\Phi^- \coprod \Phi_\mathcal{K}^+) = \Phi^- \backslash \Phi_\mathcal{I} \coprod \Phi_\mathcal{K}^+ \backslash \Phi_\mathcal{I}.$$

    Suppose $\Phi \backslash \Gamma(\mathcal{I}, \mathcal{K})$ is closed. Then there exist $\alpha \in \mathcal{K} \backslash \mathcal{I}$, $\beta \in \mathcal{I}$ such that $(\alpha, \beta) \neq 0$. Then $-\alpha - \beta \in \Phi^- \backslash \Phi_\mathcal{I}$ and $\alpha \in \Phi_\mathcal{K}^+ \backslash \Phi_\mathcal{I}$. Now $-\beta = -\alpha - \beta + \alpha \notin \Phi \backslash \Gamma(\mathcal{I}, \mathcal{K})$, which is a contradiction. Therefore, $\mathcal{I}$ is orthogonal to $\mathcal{K} \backslash \mathcal{I}$.

    If $\mathcal{I}$ is orthogonal to $\mathcal{K} \backslash \mathcal{I}$, then by Proposition \ref{general closed condition}, we can set $X = \mathcal{I}$, $Y = \mathcal{K} \backslash \mathcal{I}$, and thus $\Gamma(\mathcal{I}, \mathcal{K})$ is closed. Moreover,
    $$\Phi \backslash \Gamma(\mathcal{I}, \mathcal{K}) = \Phi^- \backslash \Phi_\mathcal{I} \coprod \Phi_\mathcal{K}^+ \backslash \Phi_\mathcal{I} = \Phi^- \backslash \Phi_\mathcal{I} \coprod \Phi^+_{\mathcal{K} \backslash \mathcal{I}}.$$

    Let $\beta \in \Phi^- \backslash \Phi_\mathcal{I}$. We can decompose $\beta$ uniquely as 
    $\beta = \beta_1 + \beta_2 + \beta_3$
    such that $\beta_1 \in \langle \Pi \backslash (\mathcal{I} \cup \mathcal{K}) \rangle$, $\beta_2 \in \langle \mathcal{K} \backslash \mathcal{I} \rangle$, $\beta_3 \in \langle \mathcal{I} \rangle$. Let $\alpha \in \Phi^+_{\mathcal{K} \backslash \mathcal{I}}$ and suppose $\alpha + \beta \in \Phi$. If $\beta_1 \neq 0$, then $\alpha + \beta \in \Phi^- \backslash \Phi_\mathcal{I}$. If $\beta_1 = 0$, then $\beta_2 \neq 0$ and $\beta_3 = 0$; we have $\alpha + \beta \in \Phi_{\mathcal{K} \backslash \mathcal{I}} \subset \Phi \backslash \Gamma(\mathcal{I}, \mathcal{K})$. Since $\Phi^- \backslash \Phi_\mathcal{I}$ and $\Phi^+_{\mathcal{K} \backslash \mathcal{I}}$ are both closed, $\Phi \backslash \Gamma(\mathcal{I}, \mathcal{K})$ is closed.

    In fact, if $\Phi_\mathcal{K} = \Phi_{\mathcal{I} \cap \mathcal{K}} \coprod \Phi_{\mathcal{K} \backslash \mathcal{I}}$, then $-(\Phi \backslash \Gamma(\mathcal{I}, \mathcal{K})) = \Gamma(\mathcal{K} \backslash \mathcal{I}, \mathcal{I})$.
\end{prf}

\subsection{Weyl group of the subgroup corresponding to closed subset}Let $\Gamma$ be a closed subset of $\Phi$. The Weyl group $W_{H_\Gamma} := N_{H_\Gamma}(T)/T$ of $H_\Gamma$ is a subgroup of $W$. There is an inclusion of $H_{\Gamma_s}$ into $H_\Gamma$. 

\begin{proposition}
\label{Weyl group depend on symmetric subset}
    The inclusion of $H_{\Gamma_s}$ in $H_\Gamma$ induces an isomorphism between $W_{H_{\Gamma_s}}$ and $W_{H_{\Gamma}}$.
\end{proposition}

\begin{prf}
    Since $\Gamma_u$ is an ideal of $\Gamma$, $U_{\Gamma_u}$ is a normal subgroup of $H_\Gamma$. The quotient group $\bar{H} := H_{\Gamma}/U_{\Gamma_u}$ is connected with the same Lie algebra as $H_{\Gamma_s}$. Let $\pi: H_\Gamma \rightarrow H_{\Gamma}/U_{\Gamma_u}$ be the quotient homomorphism. Then $\pi$ restricted to $H_{\Gamma_s}$ is an isomorphism to $H_{\Gamma}/U_{\Gamma_u}$. The restriction of $\pi$ on $T \subset H_\Gamma$ is injective. Let $\bar{T} := \pi(T)$.

    Suppose $h$ lies in the kernel of $\pi: N_{H_{\Gamma}}(T) \rightarrow N_{\bar{H}}(\bar{T})$, then $h \in U_{\Gamma_u}$ is unipotent. Thus, $h \in N_{G}(T) \cap U_{\Gamma_u}$. We could assume $U_{\Gamma_u} \subset U_{\Phi^+} = U$. But $U \cap N_G(T) = \{1\}$, the identity of $G$. So $N_{H_{\Gamma}}(T) \rightarrow N_{\bar{H}}(\bar{T})$ is injective, and the intersection $N_{H_{\Gamma_s}}(T) \hookrightarrow N_{H_{\Gamma}}(T) \hookrightarrow N_{\bar{H}}(\bar{T})$ is an isomorphism. Thus the inclusion $N_{H_{\Gamma_s}}(T) \hookrightarrow N_{H_{\Gamma}}(T)$ is an isomorphism.
\end{prf}

Since $H_{\Gamma_s}$ is reductive, we have $W_{H_{\Gamma_s}} = W_{\Gamma_s}$.

\section{Fixed point subsets}
\label{Section Some fixed points subsets}
In this section, we construct an appropriate torus action on \( X_{e,\Gamma} :=H_\Gamma\!\backslash\!\!\backslash T^*G \underset{[\mathfrak{g}^*/G]}{\times} T^*G/\!\!/_e U_e \) and calculate the fixed point set under the action in Lemma \ref{general fixed point}. When \( \Gamma = \Gamma(\mathcal{I}, \mathcal{J}, \mathcal{K}) \), we prove the main theorem \ref{fixed point set theorem} on the combinatorial description of the fixed point set of \( X_{e,\Gamma(\mathcal{I}, \mathcal{J}, \mathcal{K})} \).

\subsection{Fixed points of homogeneous space}
We recall some standard facts (see e.g. \cite{malle2011linear}).

\begin{theorem}\label{conjugate theorem}
   Let $G$ be an algebraic group. (i) Any two Borel subgroups are conjugate. (ii) Any two maximal tori are conjugate. 
\end{theorem}

\begin{proposition}
\label{closed subset torus}
    Any closed connected subgroup of a torus is a torus. 
\end{proposition}

\begin{lemma}\label{homogeneous fixed points}
    Let $T_0$ be a connected subgroup of the maximal torus $T$ such that the centralizer of $T_0$ is $T$. Let $\Gamma$ be a closed subset. Then the fixed points of the action of $T_0$ on $G/H_\Gamma$ are given by $W_\Phi/W_{\Gamma_s}$.
\end{lemma}

\begin{prf}
    We can assume $T_0$ is closed. Let $g \in G$ and suppose $gH_\Gamma/H_\Gamma$ is a fixed point of $T_0$. Then $g^{-1}T_0g \subset H_\Gamma$. By Theorem \ref{conjugate theorem} and Proposition \ref{closed subset torus}, there exists $h \in H_\Gamma$ such that $h^{-1}g^{-1}T_0gh \subset T$. Then 
    $$h^{-1}g^{-1}Tgh = h^{-1}g^{-1}C_G(T_0)gh = C_G(h^{-1}g^{-1}T_0gh) \supset C_G(T) = T,$$
    and so $gh \in N_G(T)$. Hence, $gH_\Gamma/H_\Gamma \in N_G(T)H_\Gamma/H_\Gamma = N_G(T)/(N_G(T) \cap H_\Gamma) = W_G/W_{H_\Gamma}$. Now, using Proposition \ref{Weyl group depend on symmetric subset}, we have $W_G/W_{H_\Gamma} = W_G/W_{H_{\Gamma_s}}$.
\end{prf}

\subsection{Fixed points on generalized Slodowy variety}Let $H$ be a connected closed subgroup of $G$ containing the maximal torus $T$. Then $H = H_\Gamma$ for some closed subset $\Gamma \subset \Phi$. 

In this subsection, we consider the case 
$$H \backslash\!\!\backslash T^*G \underset{[\mathfrak{g}^*/G]}{\times}  T^*G / \!\!/_{e} U_{e} = T^*(G/H) \times_{\mathfrak{g}} S_e,$$
where \( e \) is a nilpotent element and \( S_e \), its associated Slodowy slice, is defined with an \( \mathfrak{sl}_2 \)-triple \( \tau = (e, h, f) \in \mathfrak{g}^{\oplus 3} \). These are Poisson transversals in \( T^*(G/H) \) in the sense of \cite{crooks2020log}, §3.1. We call them \textit{generalized Slodowy varieties}.

We have $T^*(G/H) \cong G \times_H \mathfrak{h}^\perp$. By Lemma \ref{embedding lemma}, we can view it as a subvariety 
$$X_{e,H} := \{(gH, y) \in G/H \times \mathfrak{g} \mid y \in S_e, \operatorname{Ad}_{g^{-1}} y \in \mathfrak{h}^\perp\}.$$

Let $\gamma: \operatorname{SL}_2 \rightarrow G$ be the group homomorphism induced by this $\mathfrak{sl}_2$-triple, whose differential $\mathfrak{sl}_2 \rightarrow \mathfrak{g}$ sends 
$\begin{bmatrix}
    0 & 1\\ 
    0 & 0
\end{bmatrix}$
to $e$. Restricting the group homomorphism $\gamma$ to the diagonal subgroup 
$\{\begin{bmatrix}
    t & 0\\ 
    0 & t^{-1}
\end{bmatrix} \mid t \in \mathbb{C}^\times\}$ 
we obtain a homomorphism $\gamma: \mathbb{C}^\times \rightarrow G$ such that $\operatorname{Ad}_{\gamma(t)} e = t^2 e, \forall t \in \mathbb{C}^\times$.

We can define a $\mathbb{C}^\times$-action $\varrho$ on $X_{e,H}$ such that 
$$t \cdot (gH, y) := (\gamma(t)gH, t^{-2} \cdot \operatorname{Ad}_{\gamma(t)} y), \forall t \in \mathbb{C}^\times.$$ 

The projection to $S_e$ of $X_{e,H}$ is $\mathbb{C}^\times$-equivariant, and $S_e^{\mathbb{C}^\times} = e$.

We also have a $G_e \cap G_f =: G_\tau$ action on $X_{e,H}$ by the diagonal action on each factor.

If $\tau' = (e', h', f') = \operatorname{Ad}_a((e, h, f)) = \operatorname{Ad}_a(\tau)$ for some $a \in G$, then we can identify $G_\tau$ and $G_{\tau'}$ by $\operatorname{Ad}_a$. 

\begin{proposition}
\label{conjugation of Slodowy variety}
    If $\tau' = \operatorname{Ad}_a(\tau)$ for some $a \in G$ and $H' = bHb^{-1}$ for some $b \in G$, then $X_{e,H}$ is $\mathbb{C}^\times$- and $G_\tau$-equivariantly isomorphic to $X_{e',H'}$.
\end{proposition}

\begin{prf}
    One can check that the following map is well-defined and is an isomorphism:
    $$ \begin{aligned}
        \psi_{a,b}: X_{e,H} & \rightarrow X_{\operatorname{Ad}_a(e), bHb^{-1}}\\
        (gH, y) & \mapsto (agb^{-1}(bHb^{-1}), \operatorname{Ad}_a(y)).
    \end{aligned} $$

    Let $t \in \mathbb{C}^\times$, then 
    $$ \begin{aligned}
        \psi_{a,b}(t \cdot (gH, y)) &= \psi_{a,b}((\gamma(t^{-1})gH, t^2 \cdot \operatorname{Ad}_{\gamma(t^{-1})} y))\\
        &= (a\gamma(t^{-1})gb^{-1}(bHb^{-1}), \operatorname{Ad}_a(t^2 \cdot \operatorname{Ad}_{\gamma(t^{-1})} y))\\
        &= ((a\gamma(t^{-1})a^{-1})agb^{-1}(bHb^{-1}), t^2 \cdot \operatorname{Ad}_{a\gamma(t^{-1})a^{-1}} \operatorname{Ad}_a(y))\\
        &= t \cdot \psi_{a,b}((gH, y)).
    \end{aligned} $$

    So $\psi_{a,b}$ is $\mathbb{C}^\times$-equivariant.

    Let $g' \in G_\tau$, then 
    $$ \begin{aligned}
        \psi_{a,b}((g'gH, \operatorname{Ad}_{g'}(y))) &= (ag'gb^{-1}(bHb^{-1}), \operatorname{Ad}_a \operatorname{Ad}_{g'}(y))\\
        &= (ag'a^{-1}agb^{-1}(bHb^{-1}), \operatorname{Ad}_{ag'a^{-1}} \operatorname{Ad}_a(y))\\
        &= (ag'a^{-1}) \cdot \psi_{a,b}((gH, y)).
    \end{aligned} $$

    Therefore, $\psi_{a,b}$ is $G_\tau$-equivariant.
\end{prf}

\subsubsection{Combinatorial description of fixed point sets}
In the following, we assume $e$ is a regular nilpotent element in a Levi subgroup $M$ of $H$. To study the fixed points of $X_{e,H}$, by Proposition \ref{conjugation of Slodowy variety}, we can assume $M = M_\mathcal{L}$, where $\mathcal{L} \subset \Pi$, is a standard Levi subgroup, $e = e_\mathcal{L} = \sum_{j \in \mathcal{L}} e_j$ is the sum of some of the chosen simple root vectors, and the $\mathfrak{sl}_2$-triple $\tau = (e, f, h)$ is chosen so that $f \in \mathfrak{m}$ and $h \in \mathfrak{t}$. Hence, the Lie algebra of $M_\mathcal{L}$ is equal to $\mathfrak{h}_{\Phi_\mathcal{L}}$.

Let $T_e := \{t \in T \mid \operatorname{Ad}_t e = e\} = \{t \in T \mid \alpha_j(t) = 1, j \in \mathcal{L}\}$. Then $T_e$ is a maximal torus in the reductive group $G_\tau$, and $T_e$ acts on $X_{e,H}$ as a subgroup of $G_\tau$. $T_e$ may not be connected, for example in $G=\operatorname{SL}_2$, $T_e=\{\pm I\}$ if $e$ is the regular nilpotent element.

We are now interested in the fixed points of $X_{e,H}$ under this $T_e^\circ \times \mathbb{C}^\times$-action. Since $\mathbb{C}^\times$ contracts $S_e$ to $e$, we need to first investigate the fixed points of $G/H$ acted on by $T_e^\circ \gamma(\mathbb{C}^\times)$.

\begin{lemma}
    The centralizer of $T_e^\circ \gamma(\mathbb{C}^\times)$ is $T$.
\end{lemma}
\begin{prf}
    Since $T_e^\circ$ is a torus, its centralizer is a Levi subgroup of $G$ (see \cite{malle2011linear}, Proposition 12.10). Since $\mathfrak{t}_e = \{t \in \mathfrak{t} \mid \alpha_j(t) = 0, \forall j \in \mathcal{L}\}$, we have $C_\mathfrak{g}(\mathfrak{t}_e) = \mathfrak{h}_{\Phi_\mathcal{L}}$, and $C_G(T_e^\circ) = M_\mathcal{L}$. Since $e$ is a regular nilpotent element in $\mathfrak{h}_{\Phi_\mathcal{L}}$, we have that $h$ is a regular semisimple element in $\mathfrak{h}_{\Phi_\mathcal{L}}$. Thus, the centralizer of $\gamma(\mathbb{C}^\times)$ in $M_\mathcal{L}$ is $T$.
\end{prf}

By Lemma \ref{homogeneous fixed points}, the $T_e^\circ \gamma(\mathbb{C}^\times)$-fixed points of $G/H$ is $W/W_H(T)$. Thus we have:
\begin{lemma}
\label{general fixed point}
    The fixed points of $X_{e,H}$ under the action of $T_e^\circ \times \mathbb{C}^\times$ are isomorphic to 
    $$\mathbb{X}_{\mathcal{L},\Gamma} := \{w \in W/W_{\Gamma_s} \mid w^{-1}(\alpha_j) \in \Phi \backslash (-\Gamma), j \in \mathcal{L}\}.$$
\end{lemma}
\begin{prf}
    The condition $\operatorname{Ad}_{w^{-1}}(e) \in \mathfrak{h}^\perp$ is equivalent to $w^{-1}(\alpha_j) \in \Phi \backslash (-\Gamma), j \in \mathcal{L}$. 
\end{prf}

We also denote the above set by $\mathbb{X}_{e,H}$ to emphasize that it is the fixed-point subset of $X_{e,H} = T^*G/\!\!/H \underset{[\mathfrak{g}^*/G]}{\times}  T^*G/\!\!/_{e} U_{e}$.

    \begin{remark}
    We also attempted to generalize Lemma \ref{general fixed point} to the case \( T^*M \underset{[\mathfrak{g}^*/G]}{\times}  T^*G / \!\!/_{e} U_{e} \), where \( M \) is a \( G \)-variety. In this context, \( T^*M \) is a Hamiltonian \( G \)-variety. We can define the \( \gamma(\mathbb{C}^\times) \)-action on \( T^*M \) by 
$$
t \cdot (x, \alpha_x) := (\gamma(t)x, \gamma(t^{-1})^*(t^{-2} \alpha_x)),
$$
where \( x \in M \) and \( \alpha_x \in T^*_x M \). 
However, we have not verified whether this action preserves \( T^*M \underset{[\mathfrak{g}^*/G]}{\times}  T^*G / \!\!/_{e} U_{e} \) or not. Additionally, we cannot conclude that \( M^T = M^{T_e^\circ \gamma(\mathbb{C}^\times)} \). This equality holds when the stabilizer subgroup of every point in \( M \) is of maximal rank. 
\end{remark}

\subsubsection{Parabolic Slodowy varieties}\label{Parabolic Slodowy variety}In this subsection, we consider the case 
$$P \backslash\!\!\backslash T^*G \underset{[\mathfrak{g}^*/G]}{\times}  T^*G / \!\!/_{e} U_{e} = X_{e,P},$$ 
where $P$ is a parabolic subgroup of $G$. We assume $P$ is standard, so it is associated with a subset $\mathcal{I}$ of $\Pi$, and $P = H_{P(\mathcal{I})}=P_\mathcal{I}$.

By Lemma \ref{general fixed point}, the fixed-point subset is bijective to 
$$\mathbb{X}_{\mathcal{L},P(\mathcal{I})} = \{w \in W/W_{\mathcal{I}} \mid w^{-1}(\alpha_j) \in \Phi^+ \backslash \Phi_{\mathcal{I}}, j \in \mathcal{L}\}.$$
Since $\Phi^+ \backslash \Phi_{\mathcal{I}}$ is closed, then $w^{-1}(\alpha_j) \in \Phi^+ \backslash \Phi_{\mathcal{I}}, j \in \mathcal{L}$ is equivalent to $w^{-1}(\Phi_\mathcal{L}^+) \in \Phi^+ \backslash \Phi_{\mathcal{I}}.$ We can also naturally embed $\mathbb{X}_{\mathcal{L},P(\mathcal{I})}$ into the parabolic double coset space $W_\mathcal{L} \backslash W / W_\mathcal{I}$.

\begin{lemma}\label{coset to double coset}
    The following quotient map is bijective 
    $$\mathbb{X}_{\mathcal{L},P(\mathcal{I})} = \{w \in W/W_{\mathcal{I}} \mid w^{-1}(\Phi_\mathcal{L}^+) \subset \Phi^+ \backslash \Phi_{\mathcal{I}}\} \rightarrow \{w \in W_\mathcal{L} \backslash W / W_\mathcal{I} \mid w^{-1}(\Phi_\mathcal{L}) \subset \Phi \backslash \Phi_\mathcal{I}\}.$$
\end{lemma}
\begin{prf}
    Let $w \in W_\mathcal{L} w W_\mathcal{I}$ such that $w^{-1}(\Phi_\mathcal{L}) \subset \Phi \backslash \Phi_\mathcal{I}$. Then $w(\Phi^+ \backslash \Phi_\mathcal{I}) \cap \Phi_\mathcal{L}$ is a closed subset of $\Phi_\mathcal{L}$, and $(w(\Phi^+ \backslash \Phi_\mathcal{I}) \cap \Phi_\mathcal{L}, - (w(\Phi^+ \backslash \Phi_\mathcal{I}) \cap \Phi_\mathcal{L}))$ is a partition of $\Phi_\mathcal{L}$. By Chapter VI.1.7, Corollary 1 of \cite{bourbaki1994lie}, $w(\Phi^+ \backslash \Phi_\mathcal{I}) \cap \Phi_\mathcal{L}$ is exactly the set of positive root vectors for a basis of $\Phi_\mathcal{L}$. So there is a unique element $w'$ in $W_\mathcal{L}$ such that $w'(w(\Phi^+ \backslash \Phi_\mathcal{I}) \cap \Phi_\mathcal{L}) = \Phi_\mathcal{L}^+$. Now $w'w \in w'w W_\mathcal{I}$ satisfies $w^{-1} w'^{-1}(\Phi_\mathcal{L}^+) \subset \Phi^+ \backslash \Phi_{\mathcal{I}}$. Therefore, the above map is one-to-one and surjective.
\end{prf}

Let $(W_\mathcal{L} \backslash W / W_\mathcal{I})^{\mathrm{free}}$ be the set of free $W_\mathcal{L} \times W_\mathcal{I}$-orbits in $W$. It is isomorphic to $(_{W_\mathcal{L} \times W_\mathcal{I}}W)/(W_\mathcal{L} \times W_\mathcal{I})$. The following lemma originally appears in the appendix of \cite{hoang2024around} and is presented without a proof. To improve clarity and completeness, we include a proof below.

\begin{lemma}[\cite{hoang2024around}]\label{double coset as free}
    $$(W_\mathcal{L} \backslash W / W_\mathcal{I})^{\mathrm{free}} = \{w \in W_\mathcal{L} \backslash W / W_\mathcal{I} \mid w^{-1}(\Phi_\mathcal{L}) \subset \Phi \backslash \Phi_\mathcal{I}\}$$
\end{lemma}
\begin{prf}
    Suppose for some $\alpha \in \Phi_\mathcal{L}$, $\beta = w^{-1}(\alpha) \in W_\mathcal{I}$. Then $s_\alpha w s_\beta = w$. And so $w$ does not lie in a free $W_\mathcal{L} \times W_\mathcal{I}$-orbit.
    
    Conversely, if $w$ does not lie in a free $W_\mathcal{L} \times W_\mathcal{I}$-orbit, then $w^{-1}W_\mathcal{L} w \cap W_\mathcal{I}$ is not a trivial subgroup of $W$. Recall that the intersection of parabolic subgroups of Coxeter groups is again a parabolic subgroup (see \cite{solomon1976mackey}, or \cite{qi2007note}), so $w^{-1}W_\mathcal{L} w \cap W_\mathcal{I} = uW_Ku^{-1},$ for some nonempty subset $K \subset \Pi$. Suppose $k \in K$, so $us_{\alpha_k}u^{-1} = s_{u(\alpha_k)} \in w^{-1}W_\mathcal{L} w \cap W_\mathcal{I}$. By Proposition 1.14 of \cite{humphreys1992reflection}, $u(\alpha_k) \in \Phi_\mathcal{I}$ and $wu(\alpha_k) \in \Phi_\mathcal{L}$, and so $w^{-1}(\Phi_\mathcal{L}) \cap \Phi_\mathcal{I} \neq \varnothing$.
\end{prf}

Combining the above two lemmas, we see there is a bijection between $\mathbb{X}_{e_\mathcal{L},P_\mathcal{I}}$ and $(W_\mathcal{L} \backslash W / W_\mathcal{I})^{\mathrm{free}}$.

\subsubsection{Combinatorial Slodowy varieties}For reasons of Lemmas \ref{coset to double coset} and \ref{double coset as free}, from now on we denote the subset $$\{w \in W/W_{\mathcal{I}} \mid w^{-1}(\Phi_\mathcal{L}^+) \in \Phi^+ \backslash \Phi_{\mathcal{I}}\}$$ 
of $W/W_\mathcal{I}$ by $(W_\mathcal{L} \backslash W/W_\mathcal{I})^{\mathrm{free}}$, and we denote the subset 
$$\{w \in W/W_{\mathcal{I}} \mid w^{-1}(\Phi_\mathcal{L}^+) \in \Phi^- \backslash \Phi_{\mathcal{I}}\}$$ 
 of $W/W_{\mathcal{I}}$ by $(\overline{W_\mathcal{L}} \backslash W/W_\mathcal{I})^{\mathrm{free}}$. This corresponds to two different ways of embedding $(W_\mathcal{L} \backslash W/W_\mathcal{I})^{\mathrm{free}}$ into $W/W_{\mathcal{I}}$.

Recall we define 
$$\Gamma(\mathcal{I}, \mathcal{J}, \mathcal{K}) = \Phi_{\mathcal{I}} \cup (\Phi^+_\mathcal{J} \backslash \Phi_\mathcal{K})$$ 
for $\mathcal{I}, \mathcal{K} \subset \mathcal{J}$. So we have 
$$
\begin{aligned}
    \Phi \backslash (-\Gamma(\mathcal{I}, \mathcal{J}, \mathcal{K})) &= \Phi \backslash (\Phi_\mathcal{I} \cup \Phi^-_\mathcal{J} \backslash \Phi_\mathcal{K})\\
    &= \Phi \backslash \Phi_\mathcal{I} \cap \Phi \backslash (\Phi^-_\mathcal{J} \cap \Phi \backslash \Phi_\mathcal{K})\\
    &= \Phi \backslash \Phi_\mathcal{I} \cap (\Phi \backslash \Phi^-_\mathcal{J} \cup \Phi_\mathcal{K})\\
    &= \Phi \backslash (\Phi_\mathcal{I} \cup \Phi^-_\mathcal{J}) \cup \Phi_\mathcal{K} \backslash \Phi_\mathcal{I}\\
    &= \Phi^+ \backslash \Phi_\mathcal{I} \coprod \Phi^- \backslash \Phi_\mathcal{J} \coprod \Phi^-_\mathcal{K} \backslash \Phi_\mathcal{I}.
\end{aligned}
$$

When $\Gamma(\mathcal{I}, \mathcal{J}, \mathcal{K})$ is closed, i.e., $\mathcal{I}, \mathcal{J}, \mathcal{K}$ satisfying the relation in Proposition \ref{general closed condition}, we call 
$$H_{\Gamma(\mathcal{I}, \mathcal{J}, \mathcal{K})} \backslash\!\!\backslash T^*G \underset{[\mathfrak{g}^*/G]}{\times}  T^*G / \!\!/_{e_\mathcal{L}} U_{e_\mathcal{L}} = X_{\mathcal{L}, \Gamma(\mathcal{I}, \mathcal{J}, \mathcal{K})}$$ 
the \textit{combinatorial Slodowy variety}. Its fixed point set is

\begin{theorem}\label{fixed point set theorem}
$$\mathbb{X}_{\mathcal{L}, \Gamma(\mathcal{I}, \mathcal{J}, \mathcal{K})} = \coprod_{\substack{\mathcal{L}_1, \mathcal{L}_2, \mathcal{L}_3 \subset \mathcal{L} \\ \mathcal{L} = \mathcal{L}_1 \coprod \mathcal{L}_2 \coprod \mathcal{L}_3}} A(\mathcal{I}, \mathcal{L}_1, \mathcal{I}) \cap \bar{A}(\mathcal{I}, \mathcal{L}_2, \mathcal{J}) \cap \bar{B}(\mathcal{I}, \mathcal{L}_3, \mathcal{K}, \mathcal{I})$$
where
$$
\begin{aligned}
    A(\mathcal{I}, \mathcal{L}_1, \mathcal{I}) &= (W_{\mathcal{L}_1} \backslash W/W_\mathcal{I})^{\mathrm{free}}\\
    \bar{A}(\mathcal{I}, \mathcal{L}_2, \mathcal{J}) &= {W/W_\mathcal{I}} \times_{W/W_\mathcal{J}} (\overline{W_{\mathcal{L}_2}} \backslash W/W_\mathcal{J})^{\mathrm{free}}\\
    \bar{B}(\mathcal{I}, \mathcal{L}_3, \mathcal{K}, \mathcal{I}) &= (\overline{W_{\mathcal{L}_3}} \backslash W/W_\mathcal{I})^{\mathrm{free}} \cap (W/W_\mathcal{I} \times_{W/W_{\mathcal{I} \cup \mathcal{K}}}((W/W_{\mathcal{Y}})^{W_{\mathcal{L}_3}}/W_\mathcal{X})
\end{aligned}$$
and $\mathcal{X},\mathcal{Y}$ satisfying $\mathcal{X}\coprod \mathcal{Y}=\mathcal{I}\cup\mathcal{K}$ are subsets defined in Proposition \ref{general closed condition}.
\end{theorem}

\begin{prf}
    $$
    \begin{aligned}
       &\mathbb{X}_{\mathcal{L}, \Gamma(\mathcal{I}, \mathcal{J}, \mathcal{K})}\\
       = &\{w \in W/W_\mathcal{I} \mid w^{-1}(\mathcal{L}) \subset \Phi^+ \backslash \Phi_\mathcal{I} \coprod \Phi^- \backslash \Phi_\mathcal{J} \coprod \Phi^-_\mathcal{K} \backslash \Phi_\mathcal{I}\}\\
       =& \coprod_{\substack{\mathcal{L}_1, \mathcal{L}_2, \mathcal{L}_3 \subset \mathcal{L} \\ \mathcal{L} = \mathcal{L}_1 \coprod \mathcal{L}_2 \coprod \mathcal{L}_3}} A(\mathcal{I}, \mathcal{L}_1, \mathcal{I}) \cap \bar{A}(\mathcal{I}, \mathcal{L}_2, \mathcal{J}) \cap \bar{B}(\mathcal{I}, \mathcal{L}_3, \mathcal{K}, \mathcal{I})
    \end{aligned}
    $$
    where
    $$
    \begin{aligned}
        A(\mathcal{I}, \mathcal{L}_1, \mathcal{I}) &= \{w \in W/W_\mathcal{I} \mid w^{-1}(\mathcal{L}_1) \subset \Phi^+ \backslash \Phi_\mathcal{I}\}\\
        \bar{A}(\mathcal{I}, \mathcal{L}_2, \mathcal{J}) &= \{w \in W/W_\mathcal{I} \mid w^{-1}(\mathcal{L}_2) \subset \Phi^- \backslash \Phi_\mathcal{J}\}\\
        \bar{B}(\mathcal{I}, \mathcal{L}_3, \mathcal{K}, \mathcal{I}) &= \{w \in W/W_\mathcal{I} \mid w^{-1}(\mathcal{L}_3) \subset \Phi^-_\mathcal{K} \backslash \Phi_\mathcal{I}\}
    \end{aligned}
    $$
    
    We have the description about the first term involving $\mathcal{L}_1$:
    $$A(\mathcal{I}, \mathcal{L}_1, \mathcal{I}) = \{w \in W/W_\mathcal{I} \mid w^{-1}(\mathcal{L}_1) \subset \Phi^+ \backslash \Phi_\mathcal{I}\} = (W_{\mathcal{L}_1} \backslash W/W_\mathcal{I})^{\mathrm{free}}.$$
    
    For the second term involving $\mathcal{L}_2$, the image after projection to $W/W_\mathcal{J}$ is $(\overline{W_{\mathcal{L}_2}} \backslash W/W_\mathcal{J})^{\mathrm{free}}$, therefore
    $$\bar{A}(\mathcal{I}, \mathcal{L}_2, \mathcal{J}) = {W/W_\mathcal{I}} \times_{W/W_\mathcal{J}} (\overline{W_{\mathcal{L}_2}} \backslash W/W_\mathcal{J})^{\mathrm{free}}.$$

    Now let $\mathcal{X},\mathcal{Y}$ be the subsets of $\mathcal{I}\cup\mathcal{J}$ defined in Proposition \ref{general closed condition}. We have $W_{\mathcal{I}\cup\mathcal{K}}=W_{\mathcal{X}\cup \mathcal{Y}}=W_\mathcal{X}W_\mathcal{Y}=W_\mathcal{Y}W_\mathcal{X}$ and $W_\mathcal{X}\subset W_\mathcal{I}, W_\mathcal{Y}\subset W_\mathcal{J}$.
    
    For the third term involving $\mathcal{L}_3$, we have $\Phi^-_\mathcal{K} \backslash \Phi_\mathcal{I} = \Phi_{\mathcal{Y} \coprod \mathcal{K} \backslash \mathcal{Y}}^- \backslash \Phi_\mathcal{I} = \Phi_\mathcal{Y}^- \backslash \Phi_\mathcal{I}$, and so 
    $$\begin{aligned}
       & \{w \in W/W_\mathcal{I} \mid w^{-1}(\mathcal{L}_3) \subset \Phi^-_\mathcal{K} \backslash \Phi_\mathcal{I}\}\\
       =& \{w \in W/W_\mathcal{I} \mid w^{-1}(\mathcal{L}_3) \subset \Phi_\mathcal{Y}\} \cap \{w \in W/W_\mathcal{I} \mid w^{-1}(\mathcal{L}_3) \subset \Phi^- \backslash \Phi_\mathcal{I}\}.
    \end{aligned}$$
The first set appearing in the intersection is 
$$
\begin{aligned}
   \{w \in W/W_\mathcal{I} \mid w^{-1}(\mathcal{L}_3) \subset \Phi_\mathcal{Y}\} &= \{w \in W/W_\mathcal{I} \mid w^{-1} W_{\mathcal{L}_3} w \subset W_\mathcal{Y}\} 
   \\   &=W/W_\mathcal{I} \times_{W/W_{\mathcal{I} \cup \mathcal{K}}} ((W/W_{\mathcal{Y}})^{W_{\mathcal{L}_3}}/W_\mathcal{X})
\end{aligned}
$$
where $(W/W_{\mathcal{Y}})^{W_{\mathcal{L}_3}}$ is the fixed point set of $W_{\mathcal{L}_3}$ acting on $W/W_{\mathcal{I} \cup \mathcal{K}}$. The quotient $(W/W_{\mathcal{Y}})^{W_{\mathcal{L}_3}}/W_\mathcal{X}$ can be viewed as a subset of $(W/W_\mathcal{Y})/W_\mathcal{X}=W/(W_\mathcal{X}W_\mathcal{Y})=W/W_{\mathcal{I}\cup \mathcal{J}}$.

And the second is
$$\{w \in W/W_\mathcal{I} \mid w^{-1}(\mathcal{L}_3) \subset \Phi^- \backslash \Phi_\mathcal{I}\} = (\overline{W_{\mathcal{L}_3}} \backslash W/W_\mathcal{I})^{\mathrm{free}}.$$

Therefore,
$$
\bar{B}(\mathcal{I}, \mathcal{L}_3, \mathcal{K}, \mathcal{I})=(\overline{W_{\mathcal{L}_3}} \backslash W/W_\mathcal{I})^{\mathrm{free}} \cap (W/W_\mathcal{I} \times_{W/W_{\mathcal{I} \cup \mathcal{K}}}((W/W_{\mathcal{Y}})^{W_{\mathcal{L}_3}}/W_\mathcal{X}).
$$
\end{prf}

One immediate application to the particular case is:

\begin{corollary}
\label{Levi fixed point set}
    If $H$ is a Levi subgroup with root system $\Phi_\mathcal{I}$, then the set of fixed points of $X_{e_\mathcal{L}, H}$ is bijective to 
    $$\mathbb{X}_{\mathcal{L}, \Phi(\mathcal{I})} = \coprod_{\substack{\mathcal{L}_1, \mathcal{L}_2 \subset \mathcal{L} \\ \mathcal{L} = \mathcal{L}_1 \coprod \mathcal{L}_2}} (W_{\mathcal{L}_1} \backslash W/W_\mathcal{I})^{\mathrm{free}} \cap (\overline{W_{\mathcal{L}_2}} \backslash W/W_\mathcal{I})^{\mathrm{free}}.$$
\end{corollary}

In particular, if $\mathcal{I} = \varnothing$, then $\mathbb{X}_{\mathcal{L}, \varnothing} = W$. So we obtain a decomposition of the Weyl group:
$$W = \coprod_{\substack{\mathcal{L}_1, \mathcal{L}_2 \subset \mathcal{L} \\ \mathcal{L} = \mathcal{L}_1 \coprod \mathcal{L}_2}} (W_{\mathcal{L}_1} \backslash W) \cap (\overline{W_{\mathcal{L}_2}} \backslash W).$$

Here the elements of $W$ correspond to the fixed points of $T^*G / \!\!/ T \underset{[\mathfrak{g}^*/G]}{\times}  T^*G / \!\!/_{-e} U_{-e}$ and $e = e_\mathcal{L}$. The coset space $W_{\mathcal{L}_1} \backslash W$ (resp. $\overline{W_{\mathcal{L}_2}} \backslash W$) should be viewed as the set of representatives of the shortest (resp. longest) element in each coset.

\section{4d mirror branes and further discussions}\label{Section 4d mirror}
Let $\mathfrak{A} = \{T^*G / \!\!/ P_\mathcal{I}, T^*G / \!\!/_{{e}_\mathcal{I}} U_{e_\mathcal{I}} \mid \mathcal{I} \subset \Pi\}$ and $\langlands{\mathfrak{A}} = \{T^*\langlands{G} / \!\!/ \langlands{P}_\mathcal{I}, T^*\langlands{G} / \!\!/_{{\langlands{e}}_\mathcal{I}} \langlands{U}_{\langlands{e}_\mathcal{I}} \mid \mathcal{I} \subset \Pi\}$. There is a bijective map $\mathfrak{F} : \mathfrak{A} \rightarrow \langlands{\mathfrak{A}}$ sending $T^*G / \!\!/ P_\mathcal{I}$ to $T^*\langlands{G} / \!\!/_{\langlands{e}_\mathcal{I}} \langlands{U}_{\langlands{e}_\mathcal{I}}$ and $T^*G / \!\!/_{{e}_\mathcal{I}} U_{e_\mathcal{I}}$ to $T^*\langlands{G} / \!\!/ \langlands{P}_\mathcal{I}$.

The sets $\mathfrak{A}$ and $\langlands{\mathfrak{A}}$ have the property that for any two varieties $X, Y$ in $\mathfrak{A}$, their intersection $X \underset{[\mathfrak{g}^*/G]}{\times}  Y$ is conjectured to be the 3d mirror to $\mathfrak{F}(X) \underset{[{\langlands{\mathfrak{g}}}^*/\langlands{G}]}{\times}  \mathfrak{F}(Y)$, as shown in Section \ref{section bow variety} and also referenced in \cite{hoang2024around}.

We aim to enlarge the sets $\mathfrak{A}$ and $\langlands{\mathfrak{A}}$. Following the discussion in Section \ref{Section intersection quiver} on the star-shaped quiver gauge theory $Q_{\bar{n}}^\star$, we conjecture that \( T^*G / \!\!/ L_\mathcal{I} \) is a 4d mirror brane to \( T^*\langlands{G} / \!\!/_{\langlands{e}_\mathcal{I}} \langlands{U} \). The quiver gauge theory \( Q^\smallstar_{\bar{n}} \) suggests that \( T^*G / \!\!/ B_{L_\mathcal{I}} \) is a 4d mirror brane to \( T^*\langlands{G} / \!\!/ U_{\langlands{P}_\mathcal{I}} \), where \( B_{L_\mathcal{I}} = H_{\Phi^+_{\mathcal{I}}} \) is the Borel subgroup of \( L_\mathcal{I} \) and \( U_{\langlands{P}_\mathcal{I}} \) is the unipotent radical of \( \langlands{P}_\mathcal{I} \).

We propose the following list of 4d mirror branes that should be continually added to:
$$
\begin{array}{c|c}
\underline{X} & \langlands{\underline{X}}\\
\hline T^*G/\!\!/P_\mathcal{I} & T^*\langlands{G}/\!\!/_{\langlands{e}_\mathcal{I}} \langlands{U}_{e_\mathcal{I}}  \\
T^*G/\!\!/L_\mathcal{I} & T^*\langlands{G}/\!\!/_{\langlands{e}_\mathcal{I}} \langlands{U}\\
T^*G/\!\!/B_{L_\mathcal{I}} & T^*\langlands{G}/\!\!/{U}_{\langlands{P}_\mathcal{I}} 
\end{array}
$$

Thus, for example, we can conjecture that \( U \backslash\!\!\backslash T^*G\underset{[\mathfrak{g}^*/G]}{\times} 
 T^*G / \!\!/ P_{\mathcal{I}} \) is a 3d mirror to \( \langlands{T} \backslash\!\!\backslash T^*\langlands{G}\underset{[{\langlands{\mathfrak{g}}}^*/\langlands{G}]}{\times}  T^*\langlands{G} / \!\!/_{\langlands{e}_\mathcal{I}} \langlands{U}_{\langlands{e}_\mathcal{I}} \). Consequently, \( U \backslash\!\!\backslash T^*G\underset{[\mathfrak{g}^*/G]}{\times} 
 T^*G / \!\!/ P_{\mathcal{I}} \) should have \( |W| \) fixed points, as noted in Corollary \ref{Levi fixed point set}. However, the set-theoretic fixed-point set is \( W / W_{\mathcal{I}} \), as indicated in Proposition \ref{NegaDimension fixed proposition}. We observe that a matching of fixed points occurs when \( \mathcal{I} = \varnothing \), meaning that \( U \backslash\!\!\backslash T^*G\underset{[\mathfrak{g}^*/G]}{\times} 
 T^*G / \!\!/ P_{\mathcal{I}} \) is a nonnegative-dimensional quotient. Therefore, in such cases, an alternative method for calculating the fixed-point set of negative-dimensional quotients should be considered. Another possible way to fix this issue is by taking affinization of $U \backslash\!\!\backslash T^*G$.

From the list, we also see that the 4d mirror branes are twisted cotangent bundles. One reasonable conjecture is that the 4d mirror brane to $T^*G / \!\!/ H_{\Gamma}$ is of the form $T^*\langlands{G} / \!\!/_{\langlands{\psi}_{\Gamma}} \langlands{U}(\Gamma)$, (which requires to take affinization when necessary,) where $H_\Gamma$ is a subgroup of maximal rank in $G$ whose roots form the closed subset $\Gamma$, and $\langlands{U}(\Gamma)$ is a unipotent subgroup of $\langlands{G}$ that depends on $\Gamma$ with a character $\langlands{\psi}_{\Gamma}$. In particular, the family of unipotent subgroups $U(\Gamma(\mathcal{I}, \mathcal{J}, \mathcal{K}))$ should satisfy that $U(\Gamma(\varnothing, \mathcal{J}, \varnothing)) = U_{P_{\mathcal{J}}}$, $U(\Gamma(\mathcal{I}, \mathcal{I}, \mathcal{K})) = U$, and $U(\Gamma(\mathcal{I}, \varnothing)) = U_{e_{\mathcal{I}}}$.




\bibliographystyle{plain}

\begin{thebibliography}{10}

\bibitem{aganagic2021elliptic}
Mina Aganagic and Andrei Okounkov.
\newblock Elliptic stable envelopes.
\newblock {\em Journal of the American Mathematical Society}, 34(1):79--133, 2021.

\bibitem{arakawa2024singularities}
Tomoyuki Arakawa, Jethro Van~Ekeren, and Anne Moreau.
\newblock Singularities of nilpotent {S}lodowy slices and collapsing levels of {W}-algebras.
\newblock In {\em Forum of Mathematics, Sigma}, volume~12, page e95. Cambridge University Press, 2024.

\bibitem{balibanu2017partial}
Ana Balibanu.
\newblock The partial compactification of the universal centralizer.
\newblock {\em arXiv preprint arXiv:1710.06327}, 2017.

\bibitem{barbasch1985unipotent}
Dan Barbasch and David~A Vogan.
\newblock Unipotent representations of complex semisimple groups.
\newblock {\em Annals of Mathematics}, 121(1):41--110, 1985.

\bibitem{ben2024relative}
David Ben-Zvi, Yiannis Sakellaridis, and Akshay Venkatesh.
\newblock Relative {L}anglands duality.
\newblock {\em arXiv preprint arXiv:2409.04677}, 2024.

\bibitem{bezrukavnikov2005equivariant}
Roman Bezrukavnikov, Michael Finkelberg, and Ivan Mirkovi{\'c}.
\newblock Equivariant homology and {K}-theory of affine {G}rassmannians and {T}oda lattices.
\newblock {\em Compositio Mathematica}, 141(3):746--768, 2005.

\bibitem{bielawski1997hyperkahler}
Roger Bielawski.
\newblock Hyperk{\"a}hler structures and group actions.
\newblock {\em Journal of the London Mathematical Society}, 55(2):400--414, 1997.

\bibitem{bielawski2017slices}
Roger Bielawski.
\newblock Slices to sums of adjoint orbits, the \text{A}tiyah-\text{Hitchin} manifold, and \text{H}ilbert schemes of points.
\newblock {\em Complex Manifolds}, 4(1):16--36, 2017.

\bibitem{bisiecki1985holomorphic}
Wojciech Bisiecki.
\newblock Holomorphic {L}agrangian bundles over flag manifolds.
\newblock {\em Journal für die reine und angewandte Mathematik}, 1985(363):174--190, 1985.

\bibitem{botta2023bow}
Tommaso~Maria Botta and Rich{\'a}rd Rim{\'a}nyi.
\newblock Bow varieties: stable envelopes and their 3d mirror symmetry.
\newblock {\em arXiv preprint arXiv:2308.07300}, 2023.

\bibitem{bourbaki1994lie}
Nicolas Bourbaki.
\newblock Lie groups and {L}ie algebras.
\newblock {\em Elements of the History of Mathematics}, pages 247--267, 1994.

\bibitem{braden2014quantizations}
Tom Braden, Anthony Licata, Nicholas Proudfoot, and Ben Webster.
\newblock {Quantizations of conical symplectic resolutions {II}: category $\mathcal{O}$ and symplectic duality}.
\newblock {\em Ast{\'e'}risque}, 384:75--179, 2016.

\bibitem{braden2012quantizations}
Tom Braden, Nicholas Proudfoot, and Ben Webster.
\newblock Quantizations of conical symplectic resolutions {I}: local and global structure.
\newblock {\em Ast{\'e}risque}, 2012.

\bibitem{braverman2016towards}
Alexander Braverman, Michael Finkelberg, and Hiraku Nakajima.
\newblock {Towards a mathematical definition of {C}oulomb branches of $3$-dimensional $\mathcal{N} = 4$ gauge theories, {II}}.
\newblock {\em Adv. Theor. Math. Phys.}, 22:1071--1147, 2018.

\bibitem{Braverman:2016pwk}
Alexander Braverman, Michael Finkelberg, and Hiraku Nakajima.
\newblock {Coulomb branches of $3d$ $\mathcal{N}=4$ quiver gauge theories and slices in the affine {G}rassmannian}.
\newblock {\em Adv. Theor. Math. Phys.}, 23:75--166, 2019.

\bibitem{bullimore2017coulomb}
Mathew Bullimore, Tudor Dimofte, and Davide Gaiotto.
\newblock The coulomb branch of 3d $\mathcal{N}= 4$ theories.
\newblock {\em Communications in Mathematical Physics}, 354:671--751, 2017.

\bibitem{bullimore2016boundaries}
Mathew Bullimore, Tudor Dimofte, Davide Gaiotto, and Justin Hilburn.
\newblock Boundaries, mirror symmetry, and symplectic duality in 3d $\mathcal{N}= 4$ gauge theory.
\newblock {\em Journal of High Energy Physics}, 2016(10):1--195, 2016.

\bibitem{cabrera2019quiver}
Santiago Cabrera, Amihay Hanany, and Rudolph Kalveks.
\newblock Quiver theories and formulae for {S}lodowy slices of classical algebras.
\newblock {\em Nuclear Physics B}, 939:308--357, 2019.

\bibitem{cabrera2017nilpotent}
Santiago Cabrera, Amihay Hanany, and Zhenghao Zhong.
\newblock Nilpotent orbits and the {C}oulomb branch of ${T}^\sigma({G})$ theories: special orthogonal vs orthogonal gauge group factors.
\newblock {\em Journal of High Energy Physics}, 2017(11):1--29, 2017.

\bibitem{chan20243d}
Ki~Fung Chan and Naichung~Conan Leung.
\newblock 3d mirror symmetry is mirror symmetry.
\newblock {\em arXiv preprint arXiv:2410.03611}, 2024.

\bibitem{chriss1997representation}
Neil Chriss and Victor Ginzburg.
\newblock {\em Representation theory and complex geometry}, volume~42.
\newblock Springer, 1997.

\bibitem{cremonesi2015t}
Stefano Cremonesi, Amihay Hanany, Noppadol Mekareeya, and Alberto Zaffaroni.
\newblock ${T_\rho^\sigma(G)}$ theories and their {H}ilbert series.
\newblock {\em Journal of High Energy Physics}, 2015(1):1--56, 2015.

\bibitem{crooks2023universal}
Peter Crooks.
\newblock Universal families of twisted cotangent bundles.
\newblock {\em arXiv preprint arXiv:2306.06439}, 2023.

\bibitem{crooks2020log}
Peter Crooks and Markus Roser.
\newblock The $\log$ symplectic geometry of {P}oisson slices.
\newblock {\em Journal of Symplectic Geometry}, 2020.

\bibitem{dancer2024complex}
Andrew Dancer, Julius~F Grimminger, Johan Martens, and Zhenghao Zhong.
\newblock Complex symplectic contractions and 3d mirrors.
\newblock {\em Journal of High Energy Physics}, 2024(11):1--24, 2024.

\bibitem{devalapurkar2024ku}
Sanath~K Devalapurkar.
\newblock Ku-theoretic spectral decompositions for spheres and projective spaces.
\newblock {\em arXiv preprint arXiv:2402.03995}, 2024.

\bibitem{dinkins20223d}
Hunter Dinkins.
\newblock 3d mirror symmetry of the cotangent bundle of the full flag variety.
\newblock {\em Letters in Mathematical Physics}, 112(5):100, 2022.

\bibitem{dokovic1994closed}
D{\v{Z}}~Dokovi{\'c}, P~Check, and J-Y H{\'e}e.
\newblock On closed subsets of root systems.
\newblock {\em Canadian Mathematical Bulletin}, 37(3):338--345, 1994.

\bibitem{fu2025affine}
Baohua Fu and Jie Liu.
\newblock The affine closure of cotangent bundles of horospherical spaces.
\newblock {\em arXiv preprint arXiv:2502.06383}, 2025.

\bibitem{Gaiotto:2008ak}
Davide Gaiotto and Edward Witten.
\newblock {{S-Duality of Boundary Conditions In N=4 Super Yang-Mills Theory}}.
\newblock {\em Adv. Theor. Math. Phys.}, 13(3):721--896, 2009.

\bibitem{gammage2022perverse}
Benjamin Gammage, Justin Hilburn, and Aaron Mazel-Gee.
\newblock Perverse schobers and 3d mirror symmetry.
\newblock {\em arXiv preprint arXiv:2202.06833}, 2022.

\bibitem{gan2002quantization}
Wee~Liang Gan and Victor Ginzburg.
\newblock Quantization of {S}lodowy slices.
\newblock {\em International Mathematics Research Notices}, 2002(5):243--255, 2002.

\bibitem{gannonproof}
Tom Gannon.
\newblock Proof of the {Ginzburg-Kazhdan} conjecture.
\newblock {\em Advances in Mathematics}, 448:109701, 2024.

\bibitem{gannon2025cotangent}
Tom Gannon.
\newblock {The Cotangent Bundle of $G/U_P$ and Kostant--Whittaker Descent}.
\newblock {\em International Mathematics Research Notices}, 2025(2):rnae285, 2025.

\bibitem{gannon2025functoriality}
Tom Gannon and Ben Webster.
\newblock Functoriality of {C}oulomb branches.
\newblock {\em arXiv preprint arXiv:2501.09962}, 2025.

\bibitem{gannon2023differential}
Tom Gannon and Harold Williams.
\newblock Differential operators on the base affine space of $ {SL}_n $ and quantized {C}oulomb branches.
\newblock {\em arXiv preprint arXiv:2312.10278}, 2023.

\bibitem{hanany2020quiver}
Amihay Hanany and Rudolph Kalveks.
\newblock Quiver theories and {H}ilbert series of classical {S}lodowy intersections.
\newblock {\em Nuclear Physics B}, 952:114939, 2020.

\bibitem{hoang2024around}
Do~Kien Hoang, Vasily Krylov, and Dmytro Matvieievskyi.
\newblock Around {Hikita-Nakajima} conjecture for nilpotent orbits and parabolic {S}lodowy varieties.
\newblock {\em arXiv preprint arXiv:2410.20512}, 2024.

\bibitem{humphreys1992reflection}
James~E Humphreys.
\newblock {\em Reflection groups and \text{C}oxeter groups}.
\newblock Number~29. Cambridge university press, 1992.

\bibitem{intriligator1996mirror}
Kenneth Intriligator and Nathan Seiberg.
\newblock Mirror symmetry in three dimensional gauge theories.
\newblock {\em Physics Letters B}, 387(3):513--519, 1996.

\bibitem{ji2023bow}
Yibo Ji.
\newblock Bow varieties as symplectic reductions of ${T}^*({G}/{P})$.
\newblock {\em arXiv preprint arXiv:2312.04696}, 2023.

\bibitem{jia2022geometry}
Boming Jia.
\newblock {\em The geometry of the affine closure of $T^*(SL_n/U)$}.
\newblock PhD thesis, The University of Chicago, 2022.

\bibitem{jin2022homological}
Xin Jin.
\newblock Homological mirror symmetry for the universal centralizers.
\newblock {\em arXiv preprint arXiv:2206.09035}, 2022.

\bibitem{kamnitzer2022symplectic}
Joel Kamnitzer.
\newblock Symplectic resolutions, symplectic duality, and {C}oulomb branches.
\newblock {\em Bulletin of the London Mathematical Society}, 54(5):1515--1551, 2022.

\bibitem{kamnitzer2019category}
Joel Kamnitzer, Peter Tingley, Ben Webster, Alex Weekes, and Oded Yacobi.
\newblock On category o for affine grassmannian slices and categorified tensor products.
\newblock {\em Proceedings of the London Mathematical Society}, 119(5):1179--1233, 2019.

\bibitem{kirillov2016quiver}
Alexander Kirillov~Jr.
\newblock {\em Quiver representations and quiver varieties}, volume 174.
\newblock American Mathematical Soc., 2016.

\bibitem{malle2011linear}
Gunter Malle and Donna Testerman.
\newblock {\em Linear algebraic groups and finite groups of Lie type}, volume 133.
\newblock Cambridge university press, 2011.

\bibitem{moore20122d}
Gregory~W. Moore and Yuji Tachikawa.
\newblock {On 2d {TQFT}s whose values are holomorphic symplectic varieties}.
\newblock {\em Proc. Symp. Pure Math.}, 85:191--208, 2012.

\bibitem{nakajima2016introduction}
Hiraku Nakajima.
\newblock Introduction to quiver varieties--for ring and representation theoriests.
\newblock {\em arXiv preprint arXiv:1611.10000}, 2016.

\bibitem{nakajima2015towards}
Hiraku Nakajima.
\newblock {Towards a mathematical definition of {C}oulomb branches of $3$-dimensional $\mathcal{N}=4$ gauge theories, {I}}.
\newblock {\em Adv. Theor. Math. Phys.}, 20:595--669, 2016.

\bibitem{nakajima2022mathematical}
Hiraku Nakajima.
\newblock A mathematical definition of {C}oulomb branches of supersymmetric gauge theories and geometric {S}atake correspondences for {K}ac-{M}oody {L}ie algebras.
\newblock {\em arXiv preprint arXiv:2201.08386}, 2022.

\bibitem{Nakajima:2024mlb}
Hiraku Nakajima.
\newblock {S-dual of Hamiltonian $\mathbf G$ spaces and relative Langlands duality}.
\newblock {\em arXiv preprint arXiv:2409.06303}, 2024.

\bibitem{nakajima2017cherkis}
Hiraku Nakajima and Yuuya Takayama.
\newblock Cherkis bow varieties and {C}oulomb branches of quiver gauge theories of affine type {A}.
\newblock {\em Selecta Mathematica}, 23:2553--2633, 2017.

\bibitem{qi2007note}
Dongwen Qi.
\newblock A note on parabolic subgroups of a {C}oxeter group.
\newblock {\em Expositiones Mathematicae}, 25(1):77--81, 2007.

\bibitem{rimanyi2020bow}
Rich{\'a}rd Rim{\'a}nyi and Yiyan Shou.
\newblock Bow varieties---geometry, combinatorics, characteristic classes.
\newblock {\em arXiv preprint arXiv:2012.07814}, 2020.

\bibitem{rimanyi2019three}
Rich{\'a}rd Rim{\'a}nyi, Andrey Smirnov, Alexander Varchenko, Zijun Zhou, et~al.
\newblock Three-dimensional mirror self-symmetry of the cotangent bundle of the full flag variety.
\newblock {\em SIGMA. Symmetry, Integrability and Geometry: Methods and Applications}, 15:093, 2019.

\bibitem{rozansky1997hyper}
Lev Rozansky and Edward Witten.
\newblock Hyper-k{\"a}hler geometry and invariants of three-manifolds.
\newblock {\em Selecta Mathematica}, 3:401--458, 1997.

\bibitem{safronov2017symplectic}
Pavel Safronov.
\newblock Symplectic implosion and the {Grothendieck-Springer} resolution.
\newblock {\em Transformation Groups}, 22:767--792, 2017.

\bibitem{solomon1976mackey}
Louis Solomon.
\newblock A \text{M}ackey formula in the group ring of a {C}oxeter group.
\newblock {\em Journal of Algebra}, 41(2):255--264, 1976.

\bibitem{strominger1996mirror}
Andrew Strominger, Shing-Tung Yau, and Eric Zaslow.
\newblock Mirror symmetry is {T}-duality.
\newblock {\em Nuclear Physics B}, 479(1-2):243--259, 1996.

\bibitem{webster2011singular}
Ben Webster.
\newblock Singular blocks of parabolic category $\mathcal{O}$ and finite $\mathrm{W}$-algebras.
\newblock {\em Journal of Pure and Applied Algebra}, 215(12):2797--2804, 2011.

\end{thebibliography}

 ~\newline

\end{document}